\documentclass[11pt]{article}
\usepackage[english]{babel}
\usepackage[utf8]{inputenc}
\usepackage{amsmath,amssymb,amsfonts,amsthm}
\usepackage{algorithmic}
\usepackage{graphicx}
\usepackage{textcomp}
\usepackage{xcolor}

\usepackage{hyperref}
\usepackage{dsfont}
\usepackage[all]{xy}

\usepackage[letterpaper,top=3cm,bottom=3cm,left=3cm,right=3cm,marginparwidth=1.75cm]{geometry}

 \newcommand{\R}{\mathbb{R}}
 
 \newcommand{\N}{\mathbb{N}}

 \newcommand{\PP}{\mathbb{P}}

\newcommand{\TT}{\mathcal{T}} 
 
 \newcommand{\ep}{\varepsilon}

\theoremstyle{plain}

\newtheorem{theorem}{Theorem}[section]
\newtheorem{lemma}[theorem]{Lemma}
\newtheorem{proposition}[theorem]{Proposition}
\newtheorem{corollary}[theorem]{Corollary}

\theoremstyle{definition}

\newtheorem{definition}[theorem]{Definition}

\theoremstyle{remark}

\newtheorem{remark}[theorem]{Remark}

\newcommand{\BT}{\begin{theorem}}
\newcommand{\ET}{\end{theorem}}

\newcommand{\BL}{\begin{lemma}}
\newcommand{\EL}{\end{lemma}}

\newcommand{\BP}{\begin{proposition}}
\newcommand{\EP}{\end{proposition}}

\newcommand{\BC}{\begin{corollary}}
\newcommand{\EC}{\end{corollary}}

\newcommand{\BD}{\begin{definition}}
\newcommand{\ED}{\end{definition}}

\newcommand{\BR}{\begin{remark}}
\newcommand{\ER}{\end{remark}}

\newcommand{\BE}{\begin{equation}}
\newcommand{\EE}{\end{equation}}
\newcommand{\la}{\lambda}

\newcommand{\exho}{\tilde{h}}

\newcommand{\RR}{{\mathbb R}}
\newcommand{\ome}{{\mathbb \omega}}

\newcommand{\ZZ}{{\mathbb Z}}
\newcommand{\NN}{{\mathbb N}}

\newcommand{\cjk}{c^{(i)}_{j,k}}

\newcommand{\pjk}{\psi^{(i)}_{j,k}} 
\newcommand{\ptjk}{\tilde{\psi}^{(i)}_{j,k}}

\newsavebox{\fmbox}
\newenvironment{fmpage}[1]
 {\begin{lrbox}{\fmbox}\begin{minipage}{#1}}
 {\end{minipage}\end{lrbox}\fbox{\usebox{\fmbox}}}

 \date{\today} 

\begin{document}

\title{Gaussian processes with random pointwise regularity}
\author{C\'eline Esser\thanks{Université de Liège, B\^at. B37 Zone Polytech, All\'ee de la D\'ecouverte 12, B-4000 Li\`ege, Belgique},  St\'ephane Jaffard\thanks{ 
        Univ Paris Est Creteil, Univ Gustave Eiffel, CNRS, LAMA UMR8050, F-94010 Creteil, France
jaffard@u-pec.fr } and B\'eatrice
Vedel \thanks{Univ Bretagne Sud, CNRS UMR 6205, LMBA, F-56000, Vannes, France}}

\maketitle

\begin{abstract}
Using wavelet series, we construct Gaussian stochastic processes whose pointwise H\"older exponent is random along a typical sample path, although at every fixed point the exponent is deterministic almost surely. We first  give a local non-sparsity condition under which the H\"older exponent  almost surely is a deterministic function. 
In the opposite direction, a lacunary construction shows that the exceptional set on which the exponent is modified can have full Hausdorff dimension. For this model we determine the marginal multifractal spectra and the bivariate spectrum of the deterministic series and its Gaussian randomization.
\end{abstract}

\section{Introduction}

The study of the  regularity of Gaussian processes has a rich history, which goes back  to the seminal works of Khintchin, Lévy and Wiener on Brownian motion  concerning its  continuity, its nowhere differentiability, and the law of the iterated logarithm. 
In the case of general Gaussian processes, the deep and extremely powerful tools which were  introduced, such as metric entropy and majorizing measures,  allowed to settle   the  issues concerning the continuity and the uniform regularity of sample paths  and  also of pointwise regularity at a given point,      yielding exact deterministic pointwise moduli of continuity,  see e.g.  the   seminal articles and the  detailed expositions  by Dudley \cite{dudley2010sample}, Fernique \cite{fernique2006regularite}, Jain, Marcus and Pisier \cite{jain1974sufficient,marcus2016random},  Ledoux \cite{ledoux1991probability},  and Talagrand \cite{talagrand2014upper} to mention but a few.  

The example supplied by  Brownian motion also motivated  the following alternative  study: How does   the pointwise regularity  of a given generic sample path evolves as a function of time?  This approach opens the possibility of   obtaining random moduli of continuity, and, in the case of non-Gaussian processes, it  leads to extremly different results as the previous approach.  In order to be more specific, we first recall the notion of pointwise regularity most commonly used.   Let $B(x_0, r)$ denote  the ball of center $x_0$ and radius $r$.

\BD\label{defi:holder_reg}
    Let $f: \RR^d \rightarrow \RR $ be a locally bounded function,  let $x_0\in\R^d$ and $\alpha\geq 0$. 
    The function $f$ belongs to $C^\alpha (x_0)$ if there exists a polynomial $P_{x_0}$ of degree less than $\alpha$ and $C,r>0$ such that
    \begin{equation*}
        \forall x \in B (x_0, r) , \quad |f(x)-P_{x_0}(x-x_0)| \leq C |x-x_0|^{\alpha}.
    \end{equation*}
    The Hölder exponent of $f$ at $x_0$ is 
    \[ h_f (x_0)=\sup\{\alpha \in\R : f\in C^{\alpha} (x_0)\}. \] 
\ED
 
In the case of a random process $X $ ($= X (x, \omega )$), two notions of pointwise regularity  can be defined.  A first possibility is to  apply the above definition sample path by sample path: For each given $\omega \in \Omega $,  one  considers the Hölder exponent of the mapping $ x  \rightarrow h_X (x, \omega )$;   we will denote it by $h_X (x)$. The alternative consists in fixing a point $x$ and consider the random variable $\omega \rightarrow h_X (x, \omega )$, which we will denote by $\exho_X (x)$.  Fubini's theorem implies that  
\[ \mbox{ for a.e. } x, \quad  \mbox{a.s.,  }  \qquad h_X (x) = \exho_X (x). \] 
This simple remark raises the following  natural questions: 
\[ { \mathcal Q}_1: \qquad \mbox{ \em For which processes do $h_X (x) $ and $ \exho_X (x)$ coincide?  } \]
The second question is less precise at this point: 
\[ { \mathcal Q}_2: \qquad \mbox{ \em  If  they do not coincide,  how much can they differ? } \]
In order to illustrate the relevance of this question, let us start by recalling two classical  examples  where the answers are well known and lead  to completely different conclusions: The first  one  is   a stable  Lévy process   of stability index $\alpha \in (0,2)$.   On one hand,  $\exho_X (x)$ is a deterministic constant function: 
\begin{equation} \label{lev}
    \forall x, \quad  \mbox{ a.s. } \qquad \exho_X (x) = 1/\alpha , 
\end{equation}
see \cite{pruitt1981growth}; 
on other hand, the sample path regularity $h_X (x)$  is a random function which takes almost everywhere the value $1/\alpha$,  but also  takes all possible values in $[0, 1/\alpha ]$  and differs from  $1/\alpha$ on  a random set of dimension 1, see \cite{Jaffard:99}.   
  This  contrasts with our second example: the  fractional Brownian motion $B_H$ of Hurst exponent $H$; this Gaussian process    satisfies  
\[  \mbox{ a.s. , } \quad  \forall x,\qquad \exho_{B_H}  (x) = h_{B_H} (x) = H .   \]

Of course, one can expect answers to questions ${ \mathcal Q}_1$ and ${ \mathcal Q}_2$ only for specific classes of processes. The setting provided by Gaussian processes is natural because  as mentioned above, $\exho_X (x)$  has been the subject of a large amount of researches and its determination can be considered as settled, so that a positive answer to 
${ \mathcal Q}_1$  would allow to transfer this collection of  results into sample path properties. An important result is that,  if $X$ is  a continuous, nowhere differentiable Gaussian field defined on $\R^d$, then, for every  $x_0\in\mathbb{R}^d$, the random variable $\exho_X(x_0)$ is deterministic, see \cite[Proposition~6.2]{ayache2007wavelet}.

 For a long time,  results concerning the regularity of stochastic processes concerned   $\exho_X (x)$. However, a  shift of interest  was  motivated by modeling issues using stochastic processes:  The experimental information  available is interpreted as  a    sample path of the process and  one expects to derive information  on the  data from  the study of the theoretical   properties of the sample paths of the process, and in particular their regularity properties.

 The sole study of $\exho $ can dramatically fail to describe the nature of sample paths. For example,   L\'evy  and related processes are commonly used   in  mathematical finance in order to model  data which display random jumps, a key feature which fails to be reflected by   \eqref{lev}.   
  Such issues  motivated the introduction of  multifractal analysis, which  was initiated by comparing regularity results of experimental turbulence data with   theoretical results concerning random cascade models proposed to model such data \cite{Mandelbrot:74}.

 Let us come back to the Brownian case;  
 Khintchine's law of the iterated logarithm  states that, at almost every point of  a sample path, Brownian motion satisfies
$$ \label{regppbrown} \forall t_0, \quad a.s.  \qquad \limsup_{ |\delta| \rightarrow 0} \frac{  | B(t_0 + \delta) -B (t_0) | }{\sqrt{ |\delta| \log ( \log (1/|\delta|) )}}  =\sqrt{2} ,$$ However  random fluctuations of   the  
modulus of continuity along  sample paths have been put in evidence:
\begin{itemize}
     \item  {\em Fast points}  where identified by Orey and Taylor \cite{orey1974often}, where $$ \label{fast}  \qquad \limsup_{ |\delta| \rightarrow 0} \frac{  | B(t_0 + \delta) -B (t_0)| }{\sqrt{ |\delta|\log (1/|\delta|)  }}  >0 . $$
      \item {\em Slow points}  were identified  by Kahane  \cite{kahane1974irregularite}  where 
$$ \label{slow}   \qquad \limsup_{ |\delta| \rightarrow 0} \frac{  | B(t_0 + \delta) -B (t_0)| }{\sqrt{ |\delta| }}  < \infty . $$
\end{itemize}
Additionally,  the Hausdorff dimensions of these exceptional random sets have been determined in \cite{orey1974often,kaufman1975large}. 
These results were later extended to   a few, closely related,  Gaussian processes such as fractional Brownian motion or  some Gaussian random  wavelet series, see \cite{EL22,khoshnevisan2007fast}.  In the opposite direction, Ayache discovered a general criterium implying that  a  stationary continuous Gaussian field  satisfies   
\[  \exists C: \quad \mbox{ a.s. } \quad  \forall x, \qquad h_X (x) = C  , \] see \cite{ayache2021monofractality}.

Let us push farther the comparison between   Brownian motion and stable Lévy processes:  
In the Brownian case, the random variations of the modulus of continuity along a sample path only are of logarithmic order; this is in sharp contrast with    L\'evy processes for which  they are, in general, of a  much larger magnitude:  the function $h_X (t)$ is  a random function of time, see  \cite{Jaffard:99}. 

The question of determining how much  $h_X$ and $\exho_X$  can differ for Gaussian processes, and in particular  whether their sample paths   can have a random Hölder exponent, was not considered for a long time,  and the first indication that this can indeed happen was obtained by Ayache in 2013: In \cite{ayache2013continuous},  he    showed that   variants of Gaussian multifractional stochastic processes  can display at least one point where their Hölder exponent is   random  with positive probability, thus showing that the answer to ${\mathcal  Q}_1 $ can be negative for some Gaussian processes. 
Our purpose in the present article is to consolidate this first result by showing that this phenomenon can happen with probability one and on a set as large as possible, i.e. on a set of Hausdorff dimension $d$ (recall that this set has to be on vanishing Lebesgue measure).

 We will explore these questions through an approach which differs from the one adopted by Ayache:  We will use  the powerful tool supplied  by   wavelet series.  Note that the use of wavelet bases for the study of stochastic processes has a long history, which  can be traced  back to the seminal  article of Ciesielski  \cite{ciesielski1961holder} which  showed the relevance of the Schauder basis for the analysis of the regularity of the sample paths of Brownian motion and other Gaussian processes, thus obtaining sharp  results in terms of Besov regularity; note that  the Schauder basis is not an orthonormal basis, but it   shares the same algorithmic (dyadic)  structure  with wavelet bases, and it follows that  the  same characterizations of (pointwise and uniform) regularity hold for regularity exponents less than one, see \cite{LM86,jaffard1996local}, and also \cite{grebenkov2016multiscale} for a recent expository paper. Additionally, regularity properties of stochastic processes derived for Brownian motion through the filtering of a pseudo-differential operator were also investigated  with the help of wavelet bases in \cite{benassi1997elliptic}.
Wavelets also became a central tool in nonparametric statistics and signal processing, where thresholding and shrinkage of empirical wavelet coefficients yield efficient adaptive estimators; see for instance \cite{donoho1995denoising,donoho1994ideal, donoho1995adapting,donoho1995wavelet}.
 
 Wavelet bases shifted from  an analysis tool to a synthesis tool, with the introduction of random wavelet series (RWS), a model where the wavelet coefficients are  all independent, and share the same distributions at each scale, see \cite{Jaffard:02,esser2025p,horst2024besov,Jaffard:00b} for the (global and pointwise) regularity properties of such processes.   Note that RWS also supply examples of stochastic processes with random pointwise exponent; however,  it is not the case when the distributions  of coefficients are Gaussian (as a consequence of Theorem \ref{theo:clust} such processes satisfy  that almost surely $h_X  $ is  a constant function of $x$.

Recently, another randomization procedure based on wavelet bases was studied in \cite{EJV24}: Starting from a deterministic wavelet series, each wavelet coefficient is multiplied by an independent standard Gaussian random variable; we refer to this procedure as Gaussian wavelet randomization (GWR), see Definition \ref{def:GWR} below. The terminology differs slightly from that of \cite{EJV24}, where such randomizations were included in the broader class of random wavelet series. The previous work mainly concerned global and uniform regularity under multiplication by unbounded i.i.d. random variables and highlights a sharp contrast with Fourier randomization: Wavelet randomization needs not have a smoothing effect and may even destroy local boundedness. Here, we focus instead on the determination of $h_X$.
The aim of the present paper is twofold. First, we show that such  Gaussian processes can display  a random  exponent $h_X$ which  almost surely differs from the deterministic  exponent $\exho_X$  on a non-empty random set. Second, we investigate how large this set can be. We show that this difference can be as dramatic as in the case of Lévy processes: the exceptional set can have Hausdorff dimension one.

Our construction is based on Gaussian wavelet randomization. Let $(\psi_{j,k})_{j,k}$ denote a wavelet basis of $L^2(\mathbb{R})$, and write
$$
f=\sum_{j,k} c_{j,k}\psi_{j,k}
$$
for the wavelet expansion of a deterministic function $f$. If $(\chi_{j,k})_{j,k}$ is a family of independent standard Gaussian random variables, we define the Gaussian wavelet randomization of $f$ by
$$
X_f=\sum_{j,k} c_{j,k}\chi_{j,k}\psi_{j,k}.
$$

As recalled above, for a Gaussian process, the pointwise H\"older exponent at a fixed deterministic point is almost surely deterministic. In the case of Gaussian wavelet randomization, we can furthermore identify this deterministic value explicitly: We prove in Proposition~\ref{prop:fixed-point} that, if $f$ is a deterministic uniformly H\"older function and $X_f$ is its Gaussian wavelet randomization, then, for every \emph{fixed} point $x\in\mathbb R$,
\begin{equation}\label{eq:regfixed}
h_{X_f}(x)=h_f(x) \quad \text{almost surely.}
\end{equation}
The exceptional null set in \eqref{eq:regfixed}, however, may depend on $x$. As we shall see, this dependence is precisely what makes it possible for the pointwise regularity to behave randomly along a typical sample path.
There is nevertheless a universal restriction on the effect of Gaussian wavelet randomization. We prove in Corollary~\ref{cor:gain-spectrum} that, almost surely,
$$
0\leq h_{X_f}(x)-h_f(x)\leq1
\qquad\text{for every }x.
$$
Thus Gaussian wavelet randomization can only increase the pointwise H\"older exponent, and in dimension one this increase is at most one. We will also show that large gains in regularity can only occur on sets of correspondingly small Hausdorff dimension.

The main result of the paper shows that, despite this restriction and despite the fixed-point identity~\eqref{eq:regfixed}, the set of points where the exponent is actually modified can be as large as possible.

\begin{theorem}\label{thm:main_intro}
There exists a deterministic lacunary wavelet series $f$ such that its Gaussian wavelet randomization $X_f$ satisfies, almost surely,
$$
\dim_{\mathcal H}\{x:h_{X_f}(x)\neq h_f(x)\}=1.
$$
\end{theorem}

The striking point in Theorem~\ref{thm:main_intro} is the opposition between the two orders of quantifiers:
$$
\forall x,\qquad
\mathbb P\bigl(h_{X_f}(x)=h_f(x)\bigr)=1,
$$
whereas
$$
\mathbb P\left(\dim_{\mathcal H}\{x:h_{X_f}(x)\neq h_f(x)\}=1 \right)=1.
$$
Thus no deterministic point sees a change of exponent with positive probability, but a typical sample path contains a random exceptional set of maximal Hausdorff dimension. A natural question is then whether this change in pointwise regularity can be detected by the usual multifractal spectrum. Recall that the multifractal spectrum of a function $g$ is defined by
$$
D_g(h):=\dim_{\mathcal H}\{x: h_g(x)=h\}.
$$
It measures the Hausdorff dimension of the set of points where $g$ has pointwise H\"older exponent $H$.
A remarkable feature of our example is that the answer is negative: the randomization effect is completely invisible at the level of the usual multifractal spectra. Indeed, as a consequence of  Theorem~\ref{thm:bivariate-spectrum} we have
$$
D_f(h)=D_{X_f}(h)
$$
for every $h$.

It is therefore natural to introduce quantities which measure directly the effect of the randomization on the pointwise H\"older exponent. We first consider the \emph{randomization effect}
$$
\Delta_f(x):=h_{X_f}(x)-h_f(x)
$$
and its spectrum
$$
D_\Delta(H):=\dim_{\mathcal{H}}\{x:\Delta_f(x)=H\}.
$$
The spectrum $D_\Delta$ measures the Hausdorff dimension of the set of points where Gaussian wavelet randomization produces a gain $H$ in pointwise regularity.

The universal bound above is optimal. Indeed, the more precise analysis of
our lacunary construction shows that, for every $H\in(0,1)$, the parameters
can be chosen so that
$$
\dim_{\mathcal H}\{x:h_{X_f}(x)-h_f(x)=H\}>0.
$$
Hence gains arbitrarily close to the maximal value $1$ can occur on sets of positive Hausdorff dimension.

We also introduce the finer \emph{bivariate random spectrum}
$$
D_{(f,X_f)}(h_1,h_2) := \dim_{\mathcal H}\{x:h_f(x)=h_1,\ h_{X_f}(x)=h_2\}.
$$
While $D_\Delta$ only records the magnitude of the change in regularity, the bivariate spectrum keeps track of both the initial and the randomized H\"older exponents. In this sense, it contains finer information. Indeed,
$$
\{x:\Delta_f(x)=H\}
=
\bigcup_{h\geq 0}
\{x:h_f(x)=h,\ h_{X_f}(x)=h+H\}.
$$
For the lacunary model constructed below, we compute both spectra explicitly. In particular, the bivariate spectrum gives a complete description, in terms of Hausdorff dimension, of the joint behavior of the deterministic and randomized pointwise exponents, and the spectrum $D_\Delta$ of the regularity gain can be recovered from it.

The paper is organized as follows. In Section~\ref{sec_recall}, we recall the notions of pointwise regularity and wavelet analysis that will be used throughout the paper. In Section~\ref{sec_GWS}, we introduce Gaussian wavelet randomization and establish results valid for arbitrary uniformly H\"older functions, including the identification of the fixed-point exponent and general bounds on the possible gain in regularity. Section~\ref{sec:nonsparse} is devoted to the non-sparse regime, where we give conditions ensuring that Gaussian wavelet randomization preserves all pointwise H\"older exponents. In Section~\ref{sec:randomexp}, we introduce the lacunary model which produces random pointwise regularity. We first determine its deterministic regularity properties, then isolate the Gaussian selection argument on the associated tree, and finally combine both ingredients in a Cantor construction. This yields the exact bivariate spectrum and the spectrum of the randomization gain.

\section{Pointwise regularity and wavelets}\label{sec_recall}

This section recalls the few notions from pointwise regularity and wavelet analysis which are used below.

Let us start by recalling the definition of a wavelet basis \cite{Dau92,Mey90I,Mallat1998}.

 \BD  \label{defwavsyst}  Let $r \geq 0$; an  $r$-smooth orthonormal wavelet basis on $\RR^d$ is a collection of   functions  $\varphi$ and $\Psi^{(i)}$,    $ i \in \{ 1, \dots, 2^d -1 \} $,  with the following properties:\begin{itemize}
\item $\varphi$ and the $\Psi^{(i)}$ and all their  partial derivatives  up to  order $r$   have fast decay, 
\item The $\varphi (\cdot -k )$,   $k \in \ZZ^d$  together with the functions 
 $  2^{dj/2} \Psi^{(i)}  (2^j \cdot -k)$, $j \geq 0, \; k \in \ZZ^d$, $i  \in \{ 1, \dots, 2^d -1\}$,  
 form an orthonormal basis of $\RR^d$. 
\end{itemize}
  \ED

  We will use the following notations   for wavelets: 
 \[  \varphi_k = \varphi (\cdot - k), \quad  \pjk =  \Psi^{(i)}  (2^j \cdot -k) , \]

  and by $c_k$ and $c^{i}_{j,k}  $ the corresponding  wavelet  coefficients, so that 
  \BE \label{calculcoeff} c^{i}_{j,k}= 2^{dj} \int f(x) \psi^{(i)}_{j,k} (x) dx.  \EE  Note that we do not use a $L^2$ normalization for wavelets and wavelet coefficients; indeed the convention we use  leads to simpler formulations when dealing with pointwise regularity. 
 We will usually not mention the smoothness required for the wavelet basis, assuming that $r$ is ``large enough''; the precise value required can always be easily tracked, and one can also assume that the $C^\infty$ Meyer wavelets are used.  When working on the torus $\TT^d$,  {\em periodized wavelet bases} are used; they  are supplied by the union of the constant function $1_{\TT^d} $,   and the 1-periodic wavelets $\ptjk$ defined on $\RR^d$ by 
  \[  \ptjk (x) : = \sum_{l \in \ZZ^d}  \Psi^{(i)}  (2^j (x-l) -k) ,    \quad  i  \in \{ 1, \dots, 2^d -1\},  \;  j \geq 0, \; k  \cdot 2^{-j} \in  [0,1)^d . \] 

  We now consider a fixed smooth wavelet basis according to Definition  \ref{defwavsyst}. It follows that any $L^2$ function $f$ can be expanded on this basis  as 
   \begin{equation} \label{defexp}
     f(x) =    \sum_{k\in \ZZ^d}
   c_k  \; \varphi (x-k) +  \sum_{j\geq 0 } \sum_{i=1}^{2^d-1} \sum_{k\in \ZZ^d}
   \cjk \; \pjk (x),
 \end{equation} 
 where the series converges in $L^2$.  A key property of wavelet bases is that this expansion actually makes sense in many functional settings.  For example, if $f$ is a tempered distribution, and if one uses Meyer's wavelets, the wavelet expansion   converges in the sense of distributions (where  \eqref{calculcoeff} has  to be understood in the sense of the $  ({ \mathcal S},  { \mathcal S'})$  duality), see \cite{Meyer:95,stromberg2006modified} where wavelets are shown to be unconditional bases of ``most''  classical function spaces.

  In the following, a fixed smooth wavelet basis is chosen. We will use the same basis both for the deterministic function and for its Gaussian randomization.

Uniform H\"older spaces admit the following simple wavelet characterization, see \cite{Meyer:95}:
\begin{equation} \label{carac:holunif}
 f \in C^{\alpha}(\RR^d) \quad \mbox{ iff }\qquad
 \exists C>0:\quad \forall i,j,k,\qquad |\cjk|\leq C2^{-\alpha j}.
\end{equation}

A sequence of coefficients is a {\em uniform Hölder sequence } if it satisfies 
  the right hand side condition of \eqref{carac:holunif} for an $\alpha >0$. 

The  wavelet characterization of pointwise regularity can be stated more easily 
 using wavelet leaders  instead of  wavelet coefficients; their definition requires  the following notations.
    Wavelet coefficients can be indexed by   dyadic cubes: If 
\[   \lambda = \displaystyle \left[ \frac{k_1}{ 2^j},  \frac{k_1+1 }{ 2^j}\right[ \times \cdots \times \left[ \frac{k_d}{ 2^j},  \frac{k_d+1 }{ 2^j}\right[\;  ,  \quad  \mbox{ then  } \quad   c^{(i)}_\lambda = 2^{dj} \displaystyle\int_{\RR^d} f(x)  \psi^{(i)} (2^j x-k) dx . \] 
The cube  $3 \la $  is the cube homothetic to  $\la$ with the same center and three times wider, i.e.  it is the cube 
\begin{equation} \label{troilam} 3\la= \left[ \frac{k_1-1}{2^j}, \frac{k_1+2}{2^j}\right) 
\times \dots \times \left[ \frac{k_d-1}{2^j}, \frac{k_d+2}{2^j}\right).
\end{equation}

\BD 
Let   $f$ be  a     locally  bounded function.   The   {wavelet
  leaders} of $f$ are  
\BE  
\label{deflead} 
d_\lambda = \sup_{i, \lambda ' \subset 3 \lambda} | c^{(i)}_{  \lambda '} |.
\EE
\ED 

The fact that the supremum  in (\ref{deflead}) is finite is a consequence of the boundedness assumption for $f$. Let  $\lambda_j (x_0)$  denote the dyadic cube of width  $2^{-j}$  which contains $x_0$. Proposition \ref{prop:WLReg} below  is the reformulation of the general pointwise regularity criterium of \cite{jaffard1996wavelet} in the setting of wavelet leaders, see \cite{Jaffard:04b}.  Its statement requires the following notions.

\BD
A \emph{modulus of continuity} is a positive non-decreasing function $\theta:\RR^+\to\RR^+$ such that $\theta(0)=0$ and 
\[
\exists C>0:\qquad
\forall h>0,\quad \theta(2h)\leq C\theta(h).
\] 
  Let $f$ be a locally bounded function,  let $x_0\in\R^d$ and $\alpha\geq 0$. 
    The function $\theta$ is a modulus of continuity of $f$ at $x_0$ if there exists a polynomial $P_{x_0}$ of degree less than $\alpha$ and $C,r>0$ such that
    \[
        \forall x \in B(x_0, r) , \quad |f(x)-P_{x_0}(x-x_0)| \leq C \theta (|x-x_0|).\]
\ED

\noindent

\BP \label{prop:WLReg}   Let $f$ be a uniform Hölder function, and let $x_0 \in \RR^d$.  If  $\ome  $ is a pointwise modulus of continuity of $f$ at $x_0$, then 
\begin{equation} \label{hypdla}  \exists C>0 , \quad  \forall j \geq 0, \qquad  d_{\lambda_j (x_0)} \leq  C  \ome (2^{-j}) . \end{equation} 
Conversely, if \eqref{hypdla} holds, and if $\omega$ is a regular modulus of continuity, then 
$\ome (r) \log (1/r) $ is a pointwise modulus of continuity of $f$ at $x_0$. In particular, a  uniform Hölder function satisfies 
\begin{equation} \label{repun} \forall x_0, \qquad   h_f(x_0) = \liminf_{j \rightarrow + \infty}\;  \; \frac{ { \log 
 \left( d_{\lambda_j (x_0)} \right)  }}{\log (2^{-j})}. \end{equation}
\EP

 \section{Gaussian wavelet randomization} \label{sec_GWS}

We now fix a smooth wavelet basis as above and introduce the randomization procedure studied in the paper.

 \BD \label{def:GWR} Let $f$ be a tempered distribution. The Gaussian wavelet randomization (GWR) of its  wavelet series
  \eqref{defexp} is the random distribution 
   \begin{equation} \label{defrws}
     X_f(x) =   \sum_{k\in \ZZ^d}
   c_k \; \xi_k \varphi (x-k) + \sum_{j\geq 0 } \sum_{i=1}^{2^d-1} \sum_{k\in \ZZ^d}
   \cjk\;  \xi^{(i)}_{j,k} \; \pjk (x),
 \end{equation} 
 where the  $ \chi_k $ and the  $\xi^{(i)}_{j,k}$ are iid standard Gaussian random variables. 
 \ED 

\BR This procedure depends on the particular wavelet basis chosen. However, we will see that, if it is smooth enough, the regularity results that will be proven do not depend on this particular choice. In the following, we will  drop the low frequency component on the $\varphi ( x-k)$ which gives a smooth contribution, and therefore has no consequence for our study. 
 \ER

 The series \eqref{defrws} almost surely converges in the sense of tempered distributions. Indeed, by a direct Borel--Cantelli argument, the Gaussian multipliers have at most polynomial growth with respect to the wavelet indices, so that multiplication of the wavelet coefficients by these variables preserves convergence in $\mathcal S'$.

Sobolev regularity is also preserved almost surely. Indeed, using the wavelet characterization of $H^s$,
$$
\mathbb E \Bigg[\sum_{i,j,k}2^{2js}|\xi_{j,k}^{(i)}c_{j,k}^{(i)}|^2\Bigg]
=
\sum_{i,j,k}2^{2js}|c_{j,k}^{(i)}|^2<\infty.
$$

For H\"older regularity, we will use repeatedly the following elementary estimate on Gaussian maxima \cite{ciesielski1961holder}.

\begin{lemma}\label{lem:gaussian-max}
Let $(\xi_{j,k})_{j,k}$ be independent standard Gaussian random variables, with at most $C2^{dj}$ variables at generation $j\in \N$. Almost surely, there exists a constant $C>0$ such that, for all sufficiently large $j$,
\[
 \max_k |\xi_{j,k}|\leq C\sqrt{j}.
\]
\end{lemma}

The situation is different for uniform H\"older spaces. By Lemma~\ref{lem:gaussian-max} and the characterization~\eqref{carac:holunif}, the Gaussian randomization may introduce an additional factor of order $\sqrt j$ in the wavelet coefficients. Thus $X_f$ need not belong to the critical space $C^s$, even when $f\in C^s$. Nevertheless, since $\sqrt j\,2^{-\varepsilon j}$ is bounded for every $\varepsilon>0$, one has $X_f\in C^{s-\varepsilon}$ almost surely for every $\varepsilon>0$. This logarithmic loss in uniform regularity does not, however, imply a decrease of the pointwise H\"older exponent.

\subsection{Pointwise regularity}

We now turn to pointwise regularity. Combined with the characterization by wavelet leaders, Lemma~\ref{lem:gaussian-max} shows that Gaussian wavelet randomization cannot decrease pointwise H\"older regularity. More precisely, we have the following result.

\begin{proposition}\label{pointmod1}
Let $f$ be a uniformly H\"older function and let $x_0\in\mathbb T^d$.
If $\omega$ is a regular pointwise modulus of continuity of $f$ at $x_0$,
then, almost surely,
$$
\omega(r)(\log(1/r))^{3/2}
$$
is a pointwise modulus of continuity of $X_f$ at $x_0$.
Moreover, this statement holds simultaneously for every $x_0$.
\end{proposition}

\begin{proof}
Let $s>0$ be such that $f$ is uniformly $C^s$. We work on the event of probability one given by Lemma~\ref{lem:gaussian-max}. On this event, for every $x_0$, there exists $C>0$ such that
$$
\sup_{i,k}|\xi^{(i)}_{j,k}|\leq C\sqrt{j}.
$$
for all $j$ large enough. As $\omega$ is a regular pointwise modulus of continuity of $f$ at $x_0$, by Proposition~\ref{prop:WLReg}, there exists $C'>0$ such that
$$
d_{\lambda_j(x_0)} \leq C' \omega(2^{-j})
$$
for every $j$ large enough. Let usdenote by $D_\lambda$ the wavelet leaders of $X_f$. Since $f$ is $C^s$, its wavelet coefficients satisfy
$$
|c_{\lambda'}|\leq C2^{-sj'}
$$
for every $\lambda'$ of scale $j'$. On the other hand, since $\omega$ is a regular modulus of continuity, there exists $A>1$ such that, for all $j$ large enough,
$$
2^{-sAj}\leq C\omega(2^{-j}).
$$
It follows that, if $j'\geq Aj$, then
$$
|\xi_{\lambda'}c_{\lambda'}|
\leq C\sqrt{j'}\,2^{-sj'}
\leq C\omega(2^{-j}),
$$
for all $j$ large enough. It therefore remains to consider the scales $j\leq j'\leq Aj$. For such scales, if $\lambda'\subset 3\lambda_j(x_0)$, then by the definition of the leader,
$$
|\xi_{\lambda'}c_{\lambda'}|
\leq C'' \sqrt j\,\omega(2^{-j})
$$
for some constant $C''>0$. Taking the supremum over all $\lambda'\subset 3\lambda_j(x_0)$ yields
$$
D_{\lambda_j(x_0)}
\leq C_{x_0}\sqrt j\,\omega(2^{-j}).
$$
The quantity
$$
\sqrt{j}\,\omega(2^{-j})
$$
corresponds, up to a multiplicative constant, to the modulus
$$
\widetilde\omega(r) = \omega(r)\big(\log(1/r)\big)^{1/2}.
$$
Since $\sqrt j$ corresponds to a factor $(\log(1/r))^{1/2}$ when $r=2^{-j}$, the randomized leaders are therefore bounded by
$$
\omega(r)(\log(1/r))^{1/2}.
$$
The converse part of Proposition~\ref{prop:WLReg} introduces one additional logarithmic factor. It follows that
$$
\omega(r)(\log(1/r))^{3/2}
$$
is a pointwise modulus of continuity of $X_f$ at $x_0$.
\end{proof}

\BR
The uniform regularity assumption in Proposition~\ref{pointmod1} is essential in the present setting. Without such an assumption, Gaussian wavelet randomization may even destroy local boundedness. This phenomenon, as well as the corresponding global regularity properties, was studied in~\cite{EJV24}. It also illustrates the strong contrast with Fourier randomization, for which randomization typically has a smoothing effect.
\ER

Since logarithmic corrections do not modify the pointwise H\"older exponent, we obtain direclty the following result.

\begin{corollary}\label{cor:regularity-increase}
Let $f$ be a uniformly H\"older function. Almost surely, for every $x$,
\[
h_{X_f}(x)\geq h_f(x).
\]
\end{corollary}

The previous corollary gives an almost sure lower bound which holds simultaneously for every point. At a fixed deterministic point, one can say more: the Gaussian randomization does not modify the H\"older exponent at all.

\begin{proposition}\label{prop:fixed-point} 
Let $f$ be a uniformly H\"older function on $[0,1]^d$ and let $x_0\in[0,1]^d$ be fixed. Then
\[
 \mathbb P\bigl(h_{X_f}(x_0)=h_f(x_0)\bigr)=1.
\]
\end{proposition}

\begin{proof}
Proposition~\ref{pointmod1} gives $ h_{X_f}(x_0)\geq h_f(x_0)$ on an event of probability one. Assume first that $h_f(x_0)<+\infty$ and fix $\varepsilon>0$. By the wavelet leader characterization, there exists a sequence $j_n\to+\infty$ such that
\[
 d_{\lambda_{j_n}(x_0)}\geq 2^{-(h_f(x_0)+\varepsilon)j_n}.
\]
For each $n$, choose a dyadic interval $\lambda'_n\subset 3 \lambda_{j_n}(x_0)$  satisfying
\[
 |c_{\lambda'_n}|\geq\frac12 d_{\lambda_{j_n}(x_0)}.
\]
If $\rho>0$ is fixed, using 
\[
 \mathbb P\bigl(|Z|\leq2^{-\rho j_n}\bigr)
 \leq C2^{-\rho j_n},
\]
the Borel--Cantelli lemma implies that, almost surely,
$
 |\xi_{\lambda'_n}|\geq2^{-\rho j_n}
$
for every sufficiently large $n$. Hence the wavelet leader of $X_f$ satisfies
\[
 D_{\lambda_{j_n}(x_0)}   \geq \frac12 2^{-(h_f(x_0)+\varepsilon+\rho)j_n}
\]
for every sufficiently large $n$, hence 
$$h_{X_f}(x_0)\leq h_f(x_0)+\varepsilon+\rho.
$$
Taking the intersection over rational $\varepsilon,\rho>0$ and then letting them tend to zero yields
\[
 h_{X_f}(x_0)\leq h_f(x_0).
\]
If $h_f(x_0)=+\infty$, Proposition~\ref{pointmod1} already gives $h_{X_f}(x_0)=+\infty$. This proves the result.
\end{proof}

\begin{remark}\label{rem:quantifiers}
The order of quantifiers in Proposition~\ref{prop:fixed-point} is essential. The proposition states
\[
 \forall x_0,\qquad
 \mathbb P\bigl(h_{X_f}(x_0)=h_f(x_0)\bigr)=1,
\]
but it does not assert that there is one event of probability one on which the equality holds for every $x_0$. The lacunary process constructed in Section~\ref{sec:randomexp} shows that the latter statement can fail on a set of Hausdorff dimension one.
\end{remark}

\subsection{A general bivariate upper bound}\label{sec:general-bivariate}

We now give an estimate which controls how much Gaussian randomization can increase a pointwise exponent. For $h\geq0$, define the cumulative leader profile
\begin{equation}\label{eq:def-nuf}
 \widetilde\nu_f(h)
 :=\lim_{\varepsilon\downarrow0}\limsup_{j\to\infty}
 \frac{1}{j}\log_2
 \#\left\{\lambda\in\Lambda_j:
 d_\lambda\geq2^{-(h+\varepsilon)j}\right\}.
\end{equation}
where  $\Lambda_j$ denote the collection of dyadic cubes of generation $j$ (i.e. of width $2^{-j}$).  This quantity is closely related to the large deviation spectrum based on wavelet leaders, see e.g.~\cite{Bastin:16, Esser:17}. By construction, $\widetilde\nu_f$ is non-decreasing and right-continuous, and it takes its values in $\{-\infty\}\cup[0,1]$.

Let us also recall that the cumulative leader profile gives an upper bound for the multifractal spectrum: for every uniformly H\"older function $g$,
$$
D_g(h)\leq\widetilde\nu_g(h).
$$

The increase of pointwise regularity can only occur when some of the Gaussian variables corresponding to large wavelet coefficients are sufficiently small. We will use the following elementary estimates for Gaussian variables near the origin.

\begin{lemma}\label{lem:estimation}
For all $h_2>h_1>0$ and every $\varepsilon>0$, there exists $J\in\N$ such that, for all $j\geq J$,
\[
(1-\varepsilon)\sqrt{\frac{2}{\pi}}\,2^{-h_1j}\leq \PP(2^{-h_2 j}\leq |Z|\leq 2^{-h_1j})\leq (1+\varepsilon)\sqrt{\frac{2}{\pi}}\,2^{-h_1j},
\]
where $Z\sim\mathcal N(0,1)$.
\end{lemma}
\begin{proof}
By symmetry of the Gaussian density,
\[
\PP(2^{-h_2 j}\leq |Z|\leq 2^{-h_1j})=\sqrt{\frac{2}{\pi}}\int_{2^{-h_2j}}^{2^{-h_1j}}e^{-x^2/2}\,dx.
\]
Let $\varepsilon>0$. For $j$ large enough, one has $1-\varepsilon\leq e^{-x^2/2}\leq 1$ for every $x\in[0,2^{-h_1j}]$. Hence
\[
(1-\varepsilon)\sqrt{\frac{2}{\pi}}\left(2^{-h_1j}-2^{-h_2j}\right)\leq \PP(2^{-h_2 j}\leq |Z|\leq 2^{-h_1j})\leq \sqrt{\frac{2}{\pi}}\left(2^{-h_1j}-2^{-h_2j}\right).
\]
Since $h_2>h_1$, for $j$ large enough $2^{-h_2j}\leq \varepsilon 2^{-h_1j}$. The conclusion follows.
\end{proof}

\BP\label{prop:general-bivariate-upper}
Almost surely, one has
\begin{itemize}
    \item for every $x$,
$$
0\leq h_{X_f}(x)-h_f(x)\leq1;
$$
    \item if $h_2>h_1$,
$$
\dim_{\mathcal H}
\big\{x:h_f(x)\leq h_1,\ h_{X_f}(x)\geq h_2\big\}
\leq
\widetilde\nu_f(h_1)+h_1-h_2
\leq
1+h_1-h_2,
$$
and the set is empty if
$$
h_2-h_1>\widetilde\nu_f(h_1).
$$
\end{itemize}
\EP

\begin{proof}
The inequality $h_{X_f}(x)\geq h_f(x)$ follows from Corollary~\ref{cor:regularity-increase}. Let us prove the second point.

Fix $h_2>h_1$, $\delta>0$, and let $\varepsilon>0$ be small enough so that $h_1-h_2+4\varepsilon<0$. For every $j\in\N$, let us set
$$
F^j(h_1,h_2,\varepsilon)
=
\left\{
\lambda\in\Lambda_j:
d_\lambda\geq2^{-(h_1+\varepsilon)j}
\text{ and }
D_\lambda\leq2^{-(h_2-\varepsilon)j}
\right\},
$$
where $D_\lambda$ denotes the wavelet leaders of $X_f$. By the characterization of pointwise H\"older regularity by wavelet leaders,
$$
\{x:h_f(x)\leq h_1,\ h_{X_f}(x)\geq h_2\}
\subset
\limsup_{j\to+\infty}
\bigcup_{\lambda\in F^j(h_1,h_2,\varepsilon)}3\lambda.
$$
Hence it is enough to estimate $\#F^j(h_1,h_2,\varepsilon)$. Let
$$
A_j(h_1,\varepsilon)
=
\{\lambda\in\Lambda_j:d_\lambda\geq2^{-(h_1+\varepsilon)j}\}.
$$
By definition of $\widetilde\nu_f(h_1)$, taking $\varepsilon$ sufficiently small if necessary, one has $\#A_j(h_1,\varepsilon)\leq2^{(\widetilde\nu_f(h_1)+\delta)j}$ for $j$ large enough.

Let us fix $\lambda\in A_j(h_1,\varepsilon)$. By definition of the wavelet leader, there exists $\lambda'\subset3\lambda$ such that $|c_{\lambda'}|\geq d_\lambda/2$. Hence, for $j$ large enough,
$$
|c_{\lambda'}|>2^{-(h_1+2\varepsilon)j}.
$$
If $\lambda\in F^j(h_1,h_2,\varepsilon)$, then $|C_{\lambda'}|\leq2^{-(h_2-\varepsilon)j}$. Since $C_{\lambda'}=g_{\lambda'}c_{\lambda'}$, this implies
$
|g_{\lambda'(\lambda)}|
\leq2^{(h_1-h_2+3\varepsilon)j}.
$
Therefore, by Lemma~\ref{lem:estimation}, for $j$ large enough,
$$
\PP\left(|C_{\lambda'(\lambda)}|\leq2^{-(h_2- \varepsilon)j}\right) \leq2^{(h_1-h_2+4\varepsilon)j}.
$$
It follows that
$$
\mathbb E[\#F^j(h_1,h_2,\varepsilon)]
\leq
2^{(\widetilde\nu_f(h_1)+h_1-h_2+\delta+4\varepsilon)j}.
$$
Hence, by Markov's inequality,
$$
\PP\left(
\#F^j(h_1,h_2,\varepsilon)
\geq
2^{(\widetilde\nu_f(h_1)+h_1-h_2+\delta+5\varepsilon)j}
\right)
\leq2^{-\varepsilon j}.
$$
The Borel--Cantelli lemma implies that, almost surely, for $j$ large enough,
$$
\#F^j(h_1,h_2,\varepsilon)
\leq
2^{(\widetilde\nu_f(h_1)+h_1-h_2+\delta+5\varepsilon)j}.
$$

If $\widetilde\nu_f(h_1)+h_1-h_2+\delta+5\varepsilon<0$, then $F^j(h_1,h_2,\varepsilon)=\emptyset$ for $j$ large enough. Therefore $\{x:h_f(x)\leq h_1,\ h_{X_f}(x)\geq h_2\}$ is empty.

Otherwise, a  covering argument gives
$$
\dim_{\mathcal H}\{x:h_f(x)\leq h_1,\ h_{X_f}(x)\geq h_2\}
\leq
\widetilde\nu_f(h_1)+h_1-h_2+\delta+5\varepsilon.
$$
Letting $\varepsilon$ and $\delta$ tend to $0$ yields
$$
\dim_{\mathcal H}\{x:h_f(x)\leq h_1,\ h_{X_f}(x)\geq h_2\}
\leq
\widetilde\nu_f(h_1)+h_1-h_2.
$$

The previous estimate is obtained for fixed $h_1$ and $h_2$. Let $(h_{1,n})_n$ and $(h_{2,m})_m$ be two dense sequences. Taking the intersection of the corresponding events, it holds almost surely for every $m,n$ such that $h_{2,m}>h_{1,n}$. Now fix $h_2>h_1$ and choose $h_{1,n_k}\downarrow h_1$ and $h_{2,m_k}\uparrow h_2$. Since
$$
\{x:h_f(x)\leq h_1,\ h_{X_f}(x)\geq h_2\}
\subset
\{x:h_f(x)\leq h_{1,n_k},\ h_{X_f}(x)\geq h_{2,m_k}\},
$$
the conclusion follows by letting $k$ tend to infinity and using the right-continuity of $\widetilde\nu_f$.

Finally, since $\widetilde\nu_f(h_1)\leq1$, assume that $h_{X_f}(x)-h_f(x)>1$ for some $x$. We can choose $h_1>h_f(x)$ and $h_2<h_{X_f}(x)$ such that $h_2-h_1>1$. The previous result implies that the corresponding set is empty, which is a contradiction. Therefore $h_{X_f}(x)-h_f(x)\leq1$ for every $x$, and the first point follows.
\end{proof}

\begin{corollary}\label{cor:marginal-upper}
Almost surely, for every $h\geq0$,
\begin{equation}\label{eq:marginal-upper}
D_{X_f}(h)\leq\sup_{0\leq u\leq h}\bigl(\widetilde\nu_f(u)+u-h\bigr).
\end{equation}
\end{corollary}

\begin{proof}
By Corollary~\ref{cor:regularity-increase}, if $h_{X_f}(x)=h$, then $h_f(x)\leq h$. Therefore,
$$
\{x:h_{X_f}(x)=h\}\subset\{x:h_f(x)=h\}\cup\bigcup_{\substack{u\in\mathbb Q\\0\leq u<h}}\{x:h_f(x)\leq u,\ h_{X_f}(x)=h\}.
$$
One has $D_f(h)\leq\widetilde\nu_f(h)$.  Furthermore, for every $u<h$, Proposition~\ref{prop:general-bivariate-upper} gives
$$
\dim_{\mathcal H} \{x:h_f(x)\leq u,\ h_{X_f}(x)=h\} \leq \widetilde\nu_f(u)+u-h.
$$
Since the union is countable, we obtain
$$
D_{X_f}(h) \leq \sup_{0\leq u\leq h}\bigl(\widetilde\nu_f(u)+u-h\bigr).
$$
\end{proof}

The preceding results show that Gaussian wavelet randomization can only increase pointwise H\"older exponents, and give general restrictions on the size of this increase. In the next sections, we will study two opposite situations, according to the local abundance or sparsity of the wavelet coefficients.

\section{Pointwise regularity of non-sparse series}\label{sec:nonsparse}

We first study a condition on the wavelet coefficients which ensures that Gaussian wavelet randomization does not modify the pointwise H\"older exponents. Under this condition, almost surely, $h_{X_f}(x)=h_f(x)$ for every $x$, so that the pointwise regularity remains deterministic.

Proposition~\ref{pointmod1} shows that Gaussian wavelet randomization cannot decrease the pointwise regularity. The question is therefore to determine when it can strictly increase it. By Proposition~\ref{prop:WLReg}, this can only happen if the randomization significantly decreases the wavelet leaders along the dyadic intervals containing the point considered.

The mechanism is simple. If a wavelet leader is generated by many coefficients of comparable size, all of them have to be multiplied by small Gaussian variables in order to significantly decrease the leader. By independence, this event has very small probability. On the contrary, if the large wavelet coefficients are isolated, a single small Gaussian multiplier may be sufficient to decrease the corresponding leader. Thus, the possibility of increasing pointwise regularity depends strongly on the sparsity of the initial deterministic wavelet series.

Sparse wavelet representations, in which only a small fraction of the coefficients are significant, play a central role in signal processing and nonparametric statistics, notably in compression, denoising and thresholding methods~\cite{donoho1994ideal,donoho1995adapting, donoho1995wavelet}. Thus, the distinction between sparse and non-sparse wavelet expansions also reflects situations commonly encountered in applications.

We now give a sufficient condition which formalizes the previous heuristic. A Gaussian variable may take a very small value with a probability $p$ which is not negligible. However, by independence, the probability that this happens simultaneously for a family of $N$ Gaussian variables is of order $p^N$. We will therefore assume that, whenever a wavelet coefficient is large, sufficiently many coefficients of comparable size can be found in a small neighbourhood in the time-scale half-plane.

We now formalize this idea. Let $\lambda=\lambda(j,k)$ be a dyadic cube of scale $j$.  The {\em $\varepsilon$-neighbourhood $N_\varepsilon(\lambda)$ } of $\lambda$ is the set of dyadic cubes $\lambda'(j',k')$ such that
\[ | j-j' | \leq \varepsilon j \quad \mbox{ and }\quad  \left| \frac{k}{2^j}-\frac{k'}{2^{j'}} \right| \leq 2^{(\varepsilon-1) j}.  \]

\BD \label{def:filling} Fix a collection $(c_\lambda)_\lambda$ of  coefficients. For  $\varepsilon>0$ and $\lambda\in\Lambda_j$, define
\[
 M_\varepsilon(\lambda)  :=\#\left\{\lambda'\in N_\varepsilon(\lambda):  |c_{\lambda'}|\ge2^{-\varepsilon j}|c_\lambda|\right\}.
\]
We say that $(c_\lambda)_\lambda$ is \emph{uniformly filling} if, for every $\varepsilon>0$,
\[
\liminf_{j \to  + \infty} \inf_{\lambda\in\Lambda_j}M_\varepsilon(\lambda)\to +\infty.
\]
Equivalently, for every $\varepsilon>0$ and every $N\in\N$, there exists $J$ such that for all $j\ge J$ and all $\lambda\in\Lambda_j$ one has $M_\varepsilon(\lambda)\ge N$.
\ED

The simplest example is provided by
\[
\sum_{j,k} 2^{-Hj}\psi_{j,k}.
\]
For a suitable biorthogonal wavelet basis, its Gaussian randomization yields, up to a smooth low-frequency contribution, a fractional Brownian motion with Hurst parameter $H$; see Meyer, Sellan and Taqqu~\cite{MeyerSellanTaqqu}. More generally, wavelet-type expansions with independent Gaussian coefficients arise in several classical Gaussian models, including Brownian motion through the L\'evy--Ciesielski construction, fractional Brownian sheets~\cite{AyacheXiao2005}, and multifractional Brownian motion~\cite{AyacheTaqqu,ayache2007wavelet}. These examples show that Gaussian wavelet series  naturally occur in standard models of Gaussian processes and fields.

Typical examples where this condition does not hold are given by highly lacunary wavelet series, which, non-surprisingly,  are (essentially) the starting point of the counterexamples of Section \ref{sec:randomexp}. 

\BT \label{theo:clust} Let $f$ be a uniformly H\"older function and assume that its wavelet coefficient sequence is uniformly filling. Then, almost surely, 
\BE \label{expgene2} \quad \forall x,  \qquad  h_{X_f} (x) =  h_f(x). \EE 
\ET 

\begin{proof}
Corollary~\ref{cor:regularity-increase} gives $h_{X_f}(x)\geq h_f(x)$ simultaneously for every $x$ on one event of probability one. It remains to prove the converse inequality.

Fix $\varepsilon>0$ and $N\in\NN$. For every dyadic cube $\lambda\in\Lambda_j$, consider the event $\mathcal E_{\lambda,\varepsilon,N}$ that every Gaussian multiplier corresponding to a dyadic interval of $M_\varepsilon(\lambda)$ satisfies
$$
|\xi_{\lambda'}|\leq2^{-dj/N}j^{-2/N}.
$$
For all sufficiently large $j$, the uniformly filling condition gives at least $N$ such coefficients. Using the independence of the Gaussian variables and the estimate $\PP(|Z|\leq t)\leq\sqrt{2/\pi}\,t$, one has
$$
\PP(\mathcal E_{\lambda,\varepsilon,N})
\leq
\left(\sqrt{\frac{2}{\pi}}\,2^{-dj/N}j^{-2/N}\right)^N
=
\left(\frac{2}{\pi}\right)^{N/2}2^{-dj}j^{-2}.
$$
Since there are at most $2^{dj}$ dyadic cubes of scale $j$, it follows that
$$
\PP\left(\exists \la\in\Lambda_j:\mathcal E_{\la,\ep,N}\right)
\leq C_Nj^{-2}.
$$
The Borel--Cantelli lemma implies that almost surely, for $j$ large enough and every $\la\in\Lambda_j$, there exists $\la'\in N_\ep(\la)$ such that
$$
|c_{\la'}|\geq 2^{-\ep j}|c_\la|
\qquad\mbox{and}\qquad
|\xi_{\la'}|\geq 2^{-dj/N}j^{-2/N}.
$$
Taking the intersection of these events for every $N\in\NN$ and every $\varepsilon=1/m$, $m\in\NN$, we may assume that all these properties hold on the same event of probability one.

Let us now fix $x_0\in\mathbb [0,1]^{d}$. We only have to consider the case $h_f(x_0)<+\infty$. Let $(j_n)_n$ be such that
$$
\lim_{n\to+\infty}
\dfrac{\log d_{\la_{j_n}(x_0)}}{\log 2^{-j_n}}
=
h_f(x_0).
$$
For every $n$, choose $\la'_n\subset 3\la_{j_n}(x_0)$ such that
$
|c_{\la'_n}|\geq \frac12 d_{\la_{j_n}(x_0)},
$
and denote by $j'_n$ the scale of $\la'_n$. Since $f$ is uniformly H\"older and $h_f(x_0)<+\infty$, one has $j'_n\leq Cj_n$. We can now apply the previous result to $\la'_n$. There exists $\la''_n\in N_\ep(\la'_n)$ such that
$$
|c_{\la''_n}|\geq 2^{-\ep j'_n}|c_{\la'_n}|
\qquad\mbox{and}\qquad
|\xi_{\la''_n}|\geq 2^{-dj'_n/N}(j'_n)^{-2/N}.
$$
Therefore,
$$
|\xi_{\la''_n}c_{\la''_n}| \geq \frac12 2^{-(\ep+d/N)j'_n}(j'_n)^{-2/N} d_{\la_{j_n}(x_0)}.
$$
Using the definition of $N_\ep(\la'_n)$, the cube $\la''_n$ is contained in $3\la_{J_n}(x_0)$, where $J_n\geq \min\{j_n,(1-\ep)j'_n\}-C$. Since $j'_n\geq j_n$, we have $J_n\geq (1-\ep)j_n-C$. Hence
$$
D_{\la_{J_n}(x_0)} \geq \frac12 2^{-(\ep+d/N)j'_n}(j'_n)^{-2/N} d_{\la_{j_n}(x_0)}.
$$
Using $j'_n\leq Cj_n$ and the characterization by wavelet leaders, we obtain
$$
h_{X_f}(x_0)\leq \frac{h_f(x_0)+C(\ep+d/N)}{1-\ep}.
$$
Letting  $N\to+\infty$ and  $\ep\to0$ gives $h_{X_f}(x_0)\leq h_f(x_0)$.
\end{proof}

\BR
The theorem explains why the Gaussian randomization of a ``dense'' wavelet representation behaves very differently from the lacunary model below. The issue is not the Gaussian distribution itself but the number of comparable coefficients which are taken into account in a wavelet leader.
\ER

\section{A Gaussian process with random exponent}\label{sec:randomexp}

We now construct a deterministic wavelet series which is lacunary both in scale and in location. We work in one dimension for simplicity.  Its extension in several variables is straightforward. Let $\alpha>0$ and $0<\eta<1$, and define a sequence of scales $(j_n)_n$ by 
$$
 j_0=1,  \qquad
 j_{n+1}=\left\lceil e^{j_n}\right\rceil.
$$
At scale $j$, select approximately $2^{\eta j}$ equally spaced wavelet locations:
\begin{equation}\label{eq:Gamma}
 \Gamma_j  :=\Big\{  k_p=\big\lfloor p\,2^{(1-\eta)j}\big\rfloor: 0\le p<\left\lfloor2^{\eta j}\right\rfloor  \Big\}.
\end{equation}

Define
\begin{equation}\label{eq:f-lacunary}
 f=\sum_n\sum_{k\in\Gamma_{j_n}}2^{-\alpha j_n}\psi_{j_n,k}
\end{equation}
and its GWR
\begin{equation}\label{eq:X-lacunary}
 X_f=\sum_n\sum_{k\in\Gamma_{j_n}}
 \chi_{j_n,k}\,2^{-\alpha j_n}\psi_{j_n,k}.
\end{equation}

The lacunarity of the sequence of scales $(j_n)_n$ will be used twice. It separates the contributions of different active generations in the deterministic leader estimates, and it later makes it possible to select, inside every retained branch, many Gaussian multipliers having a prescribed size. We first describe the deterministic geometry of the model.

\subsection{The deterministic spectrum}

 The computation of the spectrum of $f$ is classical for lacunary wavelet series, see e.g.~\cite{Jaffard:02,Jaffard:00b,SaesSeuret}; we include the details in order to fix the approximation notation used later. For $\delta\in(0,1]$, put
$$
E_\delta
:=
\limsup_{n\to\infty}
\bigcup_{k\in\Gamma_{j_n}}
B(k2^{-j_n},2^{-\delta j_n}).
$$
The active points are uniformly spaced at distance of order $2^{-\eta j_n}$. Therefore, if $\delta<\eta$, the balls of radius $2^{-\delta j_n}$ cover the whole interval for all sufficiently large $n$, and hence
$$
E_\delta=[0,1]\qquad\forall\delta<\eta.
$$
For $\delta>\eta$, this is no longer the case, and the sets $E_\delta$ distinguish points according to how well they are approximated by the active points. This approximation rate determines their pointwise H\"older regularity.

The same spacing property also gives an upper bound for the pointwise regularity of $f$. Let $\ell_n=\lfloor\eta j_n\rfloor+1$. For every $\lambda\in\Lambda_{\ell_n}$, the interval $3\lambda$ contains, for $n$ large enough, an active dyadic interval $\lambda_{j_n,k}$ with $k\in\Gamma_{j_n}$. Hence
$$
d_\lambda\geq2^{-\alpha j_n}\geq2^{-\frac{\alpha}{\eta}\ell_n}.
$$
Proposition~\ref{prop:WLReg} therefore gives
$$
h_f(x)\leq\frac{\alpha}{\eta}
$$
for every $x\in[0,1]$. On the other hand, $f$ is uniformly $C^\alpha$ by \eqref{carac:holunif}, so that
$$
\alpha\leq h_f(x)\leq\frac{\alpha}{\eta}
$$
for every $x\in[0,1]$.

\BL\label{lem:Edelta}
For every $x\in[0,1]$, define
\[
 \delta(x) :=\sup\{\delta\in(0,1]:x\in E_\delta\}.
\]
Then
\[
 h_f(x)=\frac{\alpha}{\delta(x)}.
\]
\EL

\begin{proof}
Let us first assume that $x\in E_\delta$. Then there exist infinitely many $n$ and $k\in\Gamma_{j_n}$ such that
$|x-k2^{-j_n}|\leq 2^{-\delta j_n}$. Therefore, if $J_n=\lfloor\delta j_n\rfloor-1$, the dyadic interval $\lambda_{j_n,k}$ is included in $3\lambda_{J_n}(x)$. Hence
$$
d_{\lambda_{J_n}(x)}\geq |c_{\lambda_{j_n,k}}|=2^{-\alpha j_n}.
$$
Since $J_n/j_n\to\delta$, the characterization of pointwise regularity by wavelet leaders given in Proposition~\ref{prop:WLReg} yields $h_f(x)\leq\alpha/\delta$. Letting $\delta$ tend to $\delta(x)$ gives
$$
h_f(x)\leq\frac{\alpha}{\delta(x)}.
$$

Let us now prove the converse inequality. Assume that $\delta(x)<1$ and fix $\delta>\delta(x)$. Then $x\notin E_\delta$, so that, for $n$ large enough and every $k\in\Gamma_{j_n}$,
$$
|x-k2^{-j_n}|>2^{-\delta j_n}.
$$
Let $\lambda_{j_n,k}\subset 3\lambda_j(x)$ be associated with a non-zero coefficient. Since $\lambda_{j_n,k}\subset 3\lambda_j(x)$, one has
$|x-k2^{-j_n}|\leq C2^{-j}$. Therefore,
$2^{-\delta j_n}\leq C2^{-j}$, and thus
$j_n\geq (j-C')/\delta$. It follows that
$$
|c_{\lambda_{j_n,k}}|=2^{-\alpha j_n}
\leq C2^{-\alpha j/\delta}.
$$
Taking the supremum over all the non-zero coefficients appearing in the definition of the wavelet leader, we obtain
$d_{\lambda_j(x)}\leq C2^{-\alpha j/\delta}$. Proposition~\ref{prop:WLReg} then gives
$h_f(x)\geq\alpha/\delta$. Letting $\delta$ tend to $\delta(x)$ yields
$$
h_f(x)\geq\frac{\alpha}{\delta(x)}.
$$
If $\delta(x)=1$, the same lower bound follows from the uniform H\"older regularity of $f$, since $h_f(x)\geq\alpha$.
\end{proof}

Consequently, if one sets 
\[
G_{\delta} :=  \bigcap_{0<\delta' < \delta} E_{\delta'}
\setminus \bigcup_{\delta < \delta' \leq 1} E_{ \delta'} \ \text{ if }
\ \delta <1  \quad \text{ and }   \quad G_1 := \bigcap_{0<\delta' < 1} E_{\delta'}.
\]
then for $\delta \in (0,1]$, one has
$$G_{\delta} =  \left\{x \in [0,1] : \delta(x) = \delta \right\}  = \left\{x \in [0,1] : h_{f}(x) = \frac{\alpha}{\delta}\right\}.$$

In order to compute the Hausdorff dimension of the sets $G_\delta$, we will use the Mass Transference Principle obtained by Jaffard in \cite{Jaffard:00b} and generalized by \cite{Beresn:06}.

\BT[Mass transference principle]\label{thm_transference}
Let $(x_n,r_n)$ be a sequence with $x_n\in[0,1]$ and $r_n\downarrow0$, and let $0<s\le1$. If
\[
 \mathcal{L}^1\left(\limsup_{n\to\infty}B(x_n,r_n^s)\right)=1,
\]
then
\[
 \mathcal{H}^s\left(\limsup_{n\to\infty}B(x_n,r_n)\right)  =\mathcal{H}^s([0,1]).
\]
In particular the latter limsup set has Hausdorff dimension at least $s$.
\ET

\BP\label{prop:det-spectrum}
For every $\delta\in[\eta,1]$,
$$
 \dim_{\mathcal{H}}G_\delta=\frac{\eta}{\delta}.
$$
Consequently,
$$
 D_f(h)=\frac{\eta}{\alpha}h,
 \qquad h\in\left[\alpha,\frac{\alpha}{\eta}\right].
$$
\EP

\begin{proof}
At generation $j_n$, the set $E_\delta$ is covered by $\# \Gamma_{j_n} = \lfloor 2^{\eta j_n} \rfloor$ balls of radius $2^{-\delta j_n}$. Therefore, a standard covering argument yields 
$
  \dim_{\mathcal{H}} E_\delta\leq\frac{\eta}{\delta}.
$
Since $G_{\delta} \subset E_{\delta'}$ for every $\delta' < \delta$, one obtains 
\[
\dim_{\mathcal H} G_{\delta} \le \dim_{\mathcal H}(E_{\delta'}) \le \frac{\eta}{\delta'}
\]
for all $\delta' < \delta$, and the result follows by letting $\delta' \to \delta$.
On the other hand, the spacing of the active coefficients gives
$$ [0,1]\subset \limsup_{n\to\infty} \bigcup_{k\in\Gamma_{j_n}} B(k2^{-j_n},2^{-\eta j_n}).
$$
Theorem~\ref{thm_transference} yields
$$
 \mathcal H^{\eta/\delta}(E_\delta)
 =\mathcal H^{\eta/\delta}([0,1]).
$$
The computation of the lower bound of $\dim_{\mathcal{H}}G_\delta$ is then obtained with the following classical argument. We have
$$
\mathcal{H}^{\frac{\eta}{\delta}} (G_{\delta}) = \mathcal{H}^{\frac{\eta}{\delta}} \left( \bigcap_{\delta'  < \delta} E_{\delta'}  \right) - \mathcal{H}^{\frac{\eta}{\delta}} \left( \bigcup_{\delta'  > \delta} E_{\delta'} \right) =  \mathcal{H}^{\frac{\eta}{\delta}} \left( \bigcap_{\delta'  < \delta} E_{\delta'}  \right)
$$
since the $\frac{\eta}{\delta}$-dimensional Hausdorff measure of $E_{\delta'}$ vanishes
if $\delta' > \delta$ using the first part of the proof. Therefore,
$$
\mathcal{H}^{\frac{\eta}{\delta}} (G_{\delta}) = \mathcal{H}^{\frac{\eta}{\delta}} \left( \bigcap_{\delta'  < \delta} E_{\delta'} \right) \geq \mathcal{H}^{\frac{\eta}{\delta}} (E_{\delta} )>0,
$$
which implies $\dim_{\mathcal{H}}G_{\delta} \geq 
\frac{\eta}{\delta}$. 
\end{proof}

Let us end this section with the computation of the cumulative leader profile of the lacunary wavelet series $f$.

\BC\label{cor:nu_f}
For every $h\in[\alpha,\alpha/\eta]$,
\[
 \widetilde\nu_f(h)=\frac{\eta}{\alpha}h.
\]
\EC

\begin{proof}
Since $d_f \leq  \widetilde\nu_f$, it suffices to prove that $ \widetilde\nu_f(h)\leq\frac{\eta}{\alpha}h.$
Suppose that a leader $d_\lambda$ at scale $j$ satisfies
\[
 d_\lambda\geq2^{-(h+\varepsilon)j}.
\]
Then there exists a dyadic subinterval $\lambda' \subset 3 \lambda$ at a scale $j_n \geq j$ associated to a non-zero coefficient and satisfying
\[
 2^{-\alpha j_n}\geq2^{-(h+2\varepsilon)j},
\]
so
\[
 j_n\leq\frac{h+2\varepsilon}{\alpha}j.
\]
Because the active scales are lacunary, for large $j$, there is at most one relevant $j_n$ in this range. Since there are $\lfloor 2^{\eta j_n}\rfloor$ active coefficients at scale $j_n$, it follows that 
\[
 \#\{\lambda\in\Lambda_j:d_\lambda\geq2^{-(h+\varepsilon)j}\}
 \leq C2^{\eta(h+2\varepsilon)j/\alpha}.
\]
Letting $\varepsilon\downarrow 0$ proves the result.
\end{proof}

\subsection{Gaussian selection among the active coefficients}

Before constructing the sets used in the computation of the bivariate spectrum, we first isolate the probabilistic part of the argument. We select, among the active wavelet coefficients, those for which the Gaussian multiplier has a prescribed exponential size. We will show that sufficiently many such coefficients can be found inside each active dyadic interval and that they are also uniformly distributed in space.

In this section, we only use the fact that the sequence $(j_n)$ satisfies 
\begin{equation}\label{eq:lacunary-scales-selection}
\frac{j_n}{j_{n+1}}\longrightarrow0.
\end{equation}
For $\beta\geq0$, an active coefficient $\lambda\in\Gamma_{j_n}$ will be called {\em $\beta$-good} if
$$\frac12\,2^{-\beta j_n}\leq|\xi_\lambda|\leq2^{-\beta j_n}.$$

We first show that the $\beta$-good coefficients are uniformly distributed at the scale corresponding to their number.

\BL\label{lem:espaceBeta}
Let $\beta\in[0,\eta)$ and $\varepsilon>0$ be such that $\beta+\varepsilon<\eta$, and put
$$J_n(\beta,\varepsilon)=\lceil(\eta-\beta-\varepsilon)j_n\rceil.$$
Almost surely, for all sufficiently large $n$ and every dyadic interval $\lambda$ of generation $J_n(\beta,\varepsilon)$, there exists a $\beta$-good coefficient $\lambda'\in\Gamma_{j_n}$ such that $\lambda'\subset\lambda$.
\EL

\begin{proof}
Let $J_n=J_n(\beta,\varepsilon)$. Every dyadic interval $\lambda$ of generation $J_n$ contains at least $c2^{\eta j_n-J_n}$ elements of $\Gamma_{j_n}$. Since each of them is $\beta$-good with probability at least $c'2^{-\beta j_n}$, the independence of the random variables gives
$$
\PP\left(\lambda\text{ contains no $\beta$-good coefficient}\right)\leq\exp\left(-c''2^{\eta j_n-J_n-\beta j_n}\right)\leq\exp\left(-c''2^{\varepsilon j_n -1}\right).
$$
There are $2^{J_n}$ dyadic intervals of generation $J_n$, and therefore the probability that there exists an interval of scale $J_n$ that contains no $\beta$-good coefficient is bounded by 
$2^{J_n}\exp\left(-c''2^{\varepsilon j_n-1}\right).$ The conclusion follows from the Borel--Cantelli lemma.
\end{proof}

Finally, we consider the endpoint $\beta=\eta$. In this case, the expected number of $\eta$-good coefficients at scale $j_{n+1}$ contained in a given active interval of scale $j_n$ is of order $2^{-j_n}$ and therefore tends to zero. Thus, the previous argument cannot be applied directly. However, if we allow an arbitrarily small variation in the exponent, the expected number of admissible coefficients tends to infinity, which gives the following result.

\BP\label{prop:endpoint}
Let $\varepsilon\in(0,\eta)$. Almost surely, there exists $N\in\NN$ such that, for every $n\geq N$ and every $\lambda\in\Gamma_{j_n}$, there exists $\lambda'\in\Gamma_{j_{n+1}}$, $\lambda'\subset\lambda$, such that
$$2^{-(\eta+\varepsilon)j_{n+1}}\leq|\xi_{\lambda'}|\leq2^{-(\eta-\varepsilon)j_{n+1}}.$$
\EP

\begin{proof}
For every $n\geq0$, let
$$\Omega_n=\left\{\exists\lambda\in\Gamma_{j_n}:\forall\lambda'\in\Gamma_{j_{n+1}},\ \lambda'\subset\lambda,\ |\xi_{\lambda'}|\notin[2^{-(\eta+\varepsilon)j_{n+1}},2^{-(\eta-\varepsilon)j_{n+1}}]\right\}.$$
For $n$ large enough, there exists $c>0$ such that
$$\PP\left(2^{-(\eta+\varepsilon)j_{n+1}}\leq|Z|\leq2^{-(\eta-\varepsilon)j_{n+1}}\right)\geq c2^{-(\eta-\varepsilon)j_{n+1}}.$$
Moreover,
$$\#\{\lambda'\in\Gamma_{j_{n+1}}:\lambda'\subset\lambda\}\asymp2^{\eta j_{n+1}-j_n}.$$
Hence, by independence,
$$\PP(\Omega_n)\leq C2^{\eta j_n}\exp\left(-c'2^{\varepsilon j_{n+1}-j_n}\right).$$
Since $j_n/j_{n+1}\to0$, for $n$ large enough,
$$\varepsilon j_{n+1}-j_n\geq\frac{\varepsilon}{2}j_{n+1}.$$
Thus
$$\PP(\Omega_n)\leq C2^{\eta j_n}\exp\left(-c'2^{\frac{\varepsilon}{2}j_{n+1}}\right),$$
which is summable. The conclusion follows from the Borel--Cantelli lemma.
\end{proof}

\BC\label{cor:endpoint-branch}
Almost surely, there exists a sequence $(\lambda_n)_{n\geq n_0}$ such that $\lambda_n\in\Gamma_{j_n}$, $\lambda_{n+1}\subset\lambda_n$, and
$$\frac{-\log|\xi_{\lambda_n}|}{\log2^{j_n}}\longrightarrow\eta.$$
\EC

\begin{proof}
Let $(\varepsilon_m)_{m\geq1}$ be a decreasing sequence in $(0,\eta)$ converging to $0$. We apply Proposition~\ref{prop:endpoint} simultaneously to all $\varepsilon_m$. Since this is a countable family, there exists an event of probability one on which all these conclusions hold.  Let $(N_m)_m$ be an increasing sequence such that Proposition~\ref{prop:endpoint} holds with $\varepsilon=\varepsilon_m$ for every $n\geq N_m$. Starting from any $\lambda_{N_1}\in\Gamma_{j_{N_1}}$, we construct the sequence recursively. If $N_m\leq n<N_{m+1}$, we choose $\lambda_{n+1}\in\Gamma_{j_{n+1}}$, $\lambda_{n+1}\subset\lambda_n$, such that
$$2^{-(\eta+\varepsilon_m)j_{n+1}}\leq|\xi_{\lambda_{n+1}}|\leq2^{-(\eta-\varepsilon_m)j_{n+1}}.$$
Therefore,
$$\eta-\varepsilon_m\leq\frac{-\log|\xi_{\lambda_{n+1}}|}{\log 2^{j_{n+1}}}\leq\eta+\varepsilon_m$$
whenever $N_m\leq n<N_{m+1}$. Since $\varepsilon_m\to0$, the conclusion follows.
\end{proof}

\subsection{Maximal regularity of $X_f$}

\BP\label{prop:maxregularity}
Almost surely, for every $x\in[0,1]$,
$$h_{X_f}(x)\leq\frac{\alpha}{\eta}.$$
\EP

\begin{proof}
Let $\varepsilon>0$. Applying Lemma~\ref{lem:espaceBeta} with $\beta=0$, almost surely, for all sufficiently large $n$, every dyadic interval of generation
$J_n=\lceil(\eta-\varepsilon)j_n\rceil$ 
contains a coefficient $\lambda\in\Gamma_{j_n}$ such that
$|\xi_\lambda|\geq\frac12.$
Consequently, for every $x\in[0,1]$, the wavelet leader at the dyadic interval of generation $J_n$ containing $x$ satisfies
$$d_{J_n}(x)\geq \frac{1}{2}2^{-\alpha j_n}.$$
It follows from the wavelet characterization of the H\"older regularity that 
$$h_{X_f}(x)\leq\frac{\alpha}{\eta-\varepsilon}.$$
Applying this argument simultaneously to a sequence $(\varepsilon_m)_m$ decreasing to $0$ and taking the intersection of the corresponding events give the conclusion.
\end{proof}

The full spectrum of $X_f$ will be obtained as a corollary of the study of the bivariate spectrum.

\subsection{The bivariate spectrum}

We now obtain a lower bound for the bivariate spectrum by constructing Cantor-type sets on which the pointwise regularities of $f$ and $X_f$ are prescribed. We first describe the construction and compute its Hausdorff dimension, then determine the two regularities at its points.

Let $(\beta_m)_m$ be a fixed dense sequence in $[0,\eta]$ and put $\varepsilon_k=1/(k+1)$. We work on the event of probability one on which Lemma~\ref{lem:espaceBeta} holds simultaneously for every pair $(\beta_m,\varepsilon_k)$ such that $\beta_m+\varepsilon_k<\eta$, and on which Lemma~\ref{lem:gaussian-max} holds.

Fix $\beta\in[0,\eta)$ and choose a sequence $(\beta_{m_k})_k$ such that $|\beta_{m_k}-\beta|\leq\varepsilon_k$. For all sufficiently large $k$, one has $\beta_{m_k}+\varepsilon_k<\eta$. For each such $k$, Lemma~\ref{lem:espaceBeta} gives an integer $N_k$ such that its conclusion holds for the pair $(\beta_{m_k},\varepsilon_k)$ at every scale $j_n$ with $n\geq N_k$. We may assume that $(N_k)_k$ is strictly increasing. For $N_k\leq n<N_{k+1}$, put
$$J_n=\left\lceil(\eta-\beta_{m_k}-\varepsilon_k)j_n\right\rceil.$$
Since $\beta_{m_k}\to\beta$ and $\varepsilon_k\to0$, we have $J_n/j_n\to\eta-\beta$. By Lemma~\ref{lem:espaceBeta}, every dyadic interval of generation $J_n$ contains a $\beta_{m_k}$-good coefficient of $\Gamma_{j_n}$, where $k$ is determined by $N_k\leq n<N_{k+1}$. In each of these dyadic intervals, we select exactly one such coefficient and denote by $\widehat\Lambda_{j_n}(\beta)$ the resulting family.

Fix now $\delta\leq 1$ satisfying
$$\frac{\alpha+\beta}{\delta}\leq\frac{\alpha}{\eta}.$$
This condition implies $\delta\geq\eta$. We first assume that $\delta>\eta$. The case $\beta=0$ and $\delta=\eta$ will be treated separately.

We now construct the Cantor set $K(\beta,\delta)$. Since $\delta>\eta$, $J_n/j_n\to\eta-\beta>0$ and $j_{n-1}/j_n\to0$, we can choose $n_0$ large enough so that
$$
2^{-\delta j_n}\leq2^{-J_n}\qquad\text{and}\qquad L_{n-1}2^{J_n}\geq 8
$$
for every $n>n_0$, where
$$
L_{n_0}=1\qquad\text{and}\qquad L_n=\frac12\,2^{-\delta j_n}\quad\text{ for }n>n_0.
$$
We may also assume that $(J_n)_{n\geq n_0}$ is strictly increasing.

We start the construction with 
$$
\mathcal C_{n_0}=\{[0,1]\},\qquad L_{n_0}=1.
$$
Suppose that $\mathcal C_{n-1}$ has been constructed, its elements being pairwise disjoint closed intervals of common length $L_{n-1}$. Let $I\in\mathcal C_{n-1}$. We consider the dyadic intervals 
$Q_{n,k}=[k2^{-J_n},(k+1)2^{-J_n}]$
and keep as candidates those with even index  $k$ such that
$Q_{n,2k}\cup Q_{n,2k+1}\subset I.$
In each such interval $Q_{n,2k}$, we consider the coefficient of $\widehat\Lambda_{j_n}(\beta)$ selected in that interval. The number of such selected coefficients is at least
$
\frac{1}{2}L_{n-1}2^{J_n}-2.
$
We retain exactly
$$
m_n=\left\lfloor\frac{1}{4}L_{n-1}2^{J_n}\right\rfloor
$$
of them in each parent $I$. This is possible since $L_{n-1}2^{J_n}\geq 8$, and in particular $m_n\geq2$. To every retained coefficient $\lambda$, we associate the closed interval
$$I_{\lambda,n}=\left[x_\lambda+\frac12\,2^{-\delta j_n},\,x_\lambda+2^{-\delta j_n}\right]$$
of length $L_n$, where $x_\lambda$ denote the left endpoint of $\lambda$. Note that if $\lambda$ was selected in $Q_{n,2k}$, then $x_\lambda\in Q_{n,2k}$, so that
$$I_{\lambda,n}\subset Q_{n,2k}\cup Q_{n,2k+1}\subset I$$
since $2^{-\delta j_n}\leq2^{-J_n}$. 
We denote by $\mathcal C_n$ the family of all intervals obtained in this way, considering all the parent intervals of $\mathcal C_{n-1}$. Thus every parent contains exactly $m_n$ children, all of length $L_n$.

The left endpoints of two coefficients retained in distinct cells with even index are separated by at least $2^{-J_n}$. Consequently, the associated intervals are separated by at least a distance $\rho_n$, where
$$ \rho_n \geq 2^{-J_n}-L_n\geq\frac12\,2^{-J_n}.$$
We put
$$K_n(\beta,\delta)=\bigcup_{I\in\mathcal C_n}I\qquad\text{and}\qquad K(\beta,\delta)=\bigcap_{n\geq n_0}K_n(\beta,\delta).$$
By construction, $K_n(\beta,\delta)\subset K_{n-1}(\beta,\delta)$, and each $K_n(\beta,\delta)$ is a non-empty compact set. Hence $K(\beta,\delta)$ is a non-empty compact set.

\BP\label{prop:dimensionCantor}
Let $\beta \in (0, \eta)$ and $\delta\in(\eta,1]$ such that
$$
\frac{\alpha + \beta}{\delta} \leq \frac{\alpha}{\eta}.
$$
The Cantor set $K(\beta,\delta)$ constructed above satisfies
$$\dim_{\mathcal H}K(\beta,\delta)=\frac{\eta-\beta}{\delta}.$$
\EP

\begin{proof}
Let us prove the upper bound. Fix $s>(\eta-\beta)/\delta$. For every $n>n_0$, the definition of the Cantor gives the covering
$$K(\beta,\delta)\subseteq K_n(\beta,\delta)=\bigcup_{I\in\mathcal C_n}I.$$
The intervals of $\mathcal C_n$ have length $L_n=\frac12\,2^{-\delta j_n}$. Moreover, since we retain at most one coefficient in each dyadic interval of generation $J_n$, their number satisfies $\#\mathcal C_n\leq2^{J_n}$. Therefore,
$$\mathcal H^s_{L_n}\big(K(\beta,\delta)\big)\leq\sum_{I\in\mathcal C_n}|I|^s=\#\mathcal C_nL_n^s\leq2^{-s}2^{J_n-\delta s j_n}.$$
Since $L_n\to0$, we obtain $\mathcal H^s(K(\beta,\delta))=0$. Letting $s$ decrease to $(\eta-\beta)/\delta$ gives
$$\dim_{\mathcal H}K(\beta,\delta)\leq\frac{\eta-\beta}{\delta}.$$

Let us now prove the lower bound. The construction satisfies the assumptions of Falconer's estimate~\cite[Example~4.6]{Falconer:86}, which gives
$$
\dim_{\mathcal H}K(\beta,\delta)\geq\liminf_{n\to\infty}\frac{\log_2 \#\mathcal C_{n-1}}{-\log_2 (m_n\rho_n)}.
$$
By construction, $\rho_n \geq \frac{1}{2}2^{-J_n}$, while
$
m_n \geq \frac{1}{8}L_{n-1} 2^{J_n}
= \frac{1}{16} 2^{-\delta j_{n-1}}2^{J_n}.
$
Hence
$$
m_n\rho_n\geq \frac{1}{32}2^{-\delta j_{n-1}}.
$$
Moreover, we have
$$
\#\mathcal C_{n-1}=\prod_{p=n_0+1}^{n-1}m_p\geq m_{n-1}\geq \frac18 L_{n-2}2^{J_{n-1}}
=\frac1{16}2^{J_{n-1}}2^{-\delta j_{n-2}}.
$$
Consequently,
$$
\liminf_{n\to\infty}\frac{\log_2 \#\mathcal C_{n-1}}{-\log_2 (m_n\rho_n)}
\geq
\liminf_{n\to\infty}
\frac{J_{n-1}-\delta j_{n-2}-4}
{\delta j_{n-1}+5}.
$$
Since
$$
\frac{J_{n-1}}{j_{n-1}}\longrightarrow\eta-\beta
\qquad\text{and}\qquad
\frac{j_{n-2}}{j_{n-1}}\longrightarrow0,
$$
we obtain
$$
\dim_{\mathcal H}K(\beta,\delta)
\geq\frac{\eta-\beta}{\delta}.
$$
\end{proof}

We now determine the regularities of $f$ and $X_f$ on this Cantor. 

\BP\label{prop:regularityCantor}
Let $\beta\in[0,\eta)$ and let $\delta\in(\eta,1]$ such that
$$
\frac{\alpha+\beta}{\delta}\leq\frac{\alpha}{\eta}.
$$
For every $x\in K(\beta,\delta)$, one has
$$h_f(x)=\frac{\alpha}{\delta}\qquad\text{and}\qquad h_{X_f}(x)=\frac{\alpha+\beta}{\delta}.$$
\EP

\begin{proof}
Fix $x\in K(\beta,\delta)$. For every $n>n_0$, there is a unique interval of $\mathcal C_n$ containing $x$, associated with a coefficient $\lambda_n(x)\in\Gamma_{j_n}$. If $N_k\leq n<N_{k+1}$, then
$$\frac12\,2^{-\delta j_n}\leq x-x_{\lambda_n(x)}\leq2^{-\delta j_n},\qquad \frac12\,2^{-\beta_{m_k}j_n}\leq|\xi_{\lambda_n(x)}|\leq2^{-\beta_{m_k}j_n}.$$
Put $\ell_n=\lfloor\delta j_n\rfloor$. Since $\ell_n\leq j_n$ and $0\leq x-x_{\lambda_n(x)}\leq2^{-\ell_n}$, we have $\lambda_n(x)\subset3\lambda_{\ell_n}(x)$. Therefore,
$$d_{\lambda_{\ell_n}(x)}\geq2^{-\alpha j_n},\qquad D_{\lambda_{\ell_n}(x)}\geq\frac12\,2^{-(\alpha+\beta_{m_k})j_n}.$$
Since $\ell_n/j_n\to\delta$ and $\beta_{m_k}\to\beta$, Proposition~\ref{prop:WLReg} gives
$$h_f(x)\leq\frac{\alpha}{\delta}\qquad\text{and}\qquad h_{X_f}(x)\leq\frac{\alpha+\beta}{\delta}.$$

Let us prove the converse inequalities. Fix $h$ such that $0<h<(\alpha+\beta)/\delta$, and let $j$ be sufficiently large. Consider an active interval $\lambda\in\Gamma_{j_n}$ such that $\lambda\subset3\lambda_j(x)$. Then
$$|x-x_\lambda|\leq 2 \cdot 2^{-j}.$$

Assume first that $\lambda=\lambda_{j_n}(x)$. By construction, since $x\in K(\beta,\delta)$, we have
$$
|x-x_{\lambda_{j_n}(x)}|\geq\frac12\,2^{-\delta j_n},
$$
and therefore $j\leq\delta j_n+2$. Since $\beta_{m_k}\to\beta$, one has $\alpha+\beta_{m_k}\geq h\delta$ for all sufficiently large $n$ so that $N_k\leq n<N_{k+1}$. Hence
$$
|c_{\lambda_{j_n}(x)}\xi_{\lambda_{j_n}(x)}|\leq2^{-(\alpha+\beta_{m_k})j_n}\leq2^{-h\delta j_n}\leq2^{2h}2^{-hj}.
$$
Similarly,
$$
|c_{\lambda_{j_n}(x)}|=2^{-\alpha j_n}\leq2^{2\alpha/\delta}2^{-\alpha j/\delta}.
$$

Assume now that $\lambda\neq\lambda_n(x)$. The points $x_\lambda$ associated with distinct active coefficients of generation $j_n$ are separated by at least $\frac12\,2^{-\eta j_n}$. Since $\delta>\eta$ and $|x-x_{\lambda_n(x)}|\leq2^{-\delta j_n}$, we obtain
$$
|x-x_\lambda|
\geq|x_\lambda-x_{\lambda_n(x)}|-|x-x_{\lambda_n(x)}|
\geq\frac12\,2^{-\eta j_n}-2^{-\delta j_n}
\geq\frac14\,2^{-\eta j_n}.
$$
Together with $|x-x_\lambda|\leq2\cdot2^{-j}$, this gives $j\leq\eta j_n+3$. By Lemma~\ref{lem:gaussian-max},
$$
2^{hj}|c_\lambda\xi_\lambda|
\leq C2^{3h}\sqrt{j_n}\,2^{-(\alpha-h\eta)j_n}.
$$
Since
$$
h<\frac{\alpha+\beta}{\delta}\leq\frac{\alpha}{\eta},
$$
we have $\alpha-h\eta>0$, and therefore
$$
2^{hj}|c_\lambda\xi_\lambda|\leq C_h.
$$
For the deterministic coefficients, we simply have
$$
|c_\lambda|
=2^{-\alpha j_n}
\leq2^{3\alpha/\eta}2^{-\alpha j/\eta}
\leq2^{3\alpha/\eta}2^{-\alpha j/\delta}.
$$

Taking the supremum over all active intervals contained in $3\lambda_j(x)$, we obtain
$$
d_{\lambda_j(x)}\leq C2^{-\alpha j/\delta}
\qquad\text{and}\qquad
D_{\lambda_j(x)}\leq C_h2^{-hj}.
$$
Proposition~\ref{prop:WLReg} gives
$$
h_f(x)\geq\frac{\alpha}{\delta}
\qquad\text{and}\qquad
h_{X_f}(x)\geq h.
$$
Letting $h$ tend to $(\alpha+\beta)/\delta$ and combining with the upper bounds already obtained concludes the proof.
\end{proof}

Combining the dimension and regularity estimates gives the following lower bound for the bivariate spectrum.

\BC\label{cor:bivariate-lower}
Almost surely, for every $\beta\in[0,\eta)$ and every $\delta\in[\eta,1]$ such that 
$$
(\alpha+\beta)/\delta\leq\alpha/\eta, 
$$
one has
$$
\dim_{\mathcal H}\left\{x\in[0,1]:h_f(x)=\frac{\alpha}{\delta},\ h_{X_f}(x)=\frac{\alpha+\beta}{\delta}\right\}\geq\frac{\eta-\beta}{\delta}.
$$
\EC

\begin{proof}
If $\delta>\eta$, Propositions~\ref{prop:dimensionCantor} and~\ref{prop:regularityCantor} give
$$K(\beta,\delta)\subseteq\left\{x\in[0,1]:h_f(x)=\frac{\alpha}{\delta},\ h_{X_f}(x)=\frac{\alpha+\beta}{\delta}\right\},$$
with $\dim_{\mathcal H}K(\beta,\delta)=(\eta-\beta)/\delta$.

If $\delta=\eta$, the condition $(\alpha+\beta)/\delta\leq\alpha/\eta$ then implies $\beta=0$. We also work on the event of probability one on which Corollary~\ref{cor:regularity-increase} and the maximal regularity bound hold simultaneously for every $x$. On this event,
$$h_f(x)\leq h_{X_f}(x)\leq\frac{\alpha}{\eta}$$
for every $x\in[0,1]$. Consequently,
$$\left\{x:h_f(x)=\frac{\alpha}{\eta}\right\}=\left\{x:h_f(x)=h_{X_f}(x)=\frac{\alpha}{\eta}\right\}.$$
By Proposition~\ref{prop:det-spectrum}, this set has Hausdorff dimension $1=(\eta-\beta)/\delta$. 
All these arguments hold on the same event of probability one, which does not depend on $\beta$ or $\delta$. The conclusion therefore holds simultaneously for all the stated parameters.
\end{proof}

Combining Corollary~\ref{cor:bivariate-lower} with the general upper bound obtained in Proposition~\ref{prop:general-bivariate-upper}, we obtain the exact bivariate spectrum.

\BT\label{thm:bivariate-spectrum}
Almost surely, if
$$
\alpha\leq h_1\leq h_2\leq\frac{\alpha}{\eta}
\qquad\text{and}\qquad
h_2\leq\left(1+\frac{\eta}{\alpha}\right)h_1,
$$
one has 
$$
D_{(f,X_f)}(h_1,h_2)
=
\frac{\eta}{\alpha}h_1+h_1-h_2.
$$
Outside this domain, the corresponding level sets are empty.
\ET

\begin{proof}
The upper bound follows directly from Proposition~\ref{prop:general-bivariate-upper} and Corollary~\ref{cor:nu_f}. Let us prove the reverse inequality. Let $(h_1,h_2)$ satisfy the assumptions of the theorem and put
$$
\delta=\frac{\alpha}{h_1},
\qquad
\beta=\alpha\left(\frac{h_2}{h_1}-1\right).
$$
Then $\delta\in[\eta,1]$, $\beta\in[0,\eta]$ and
$$
\frac{\alpha+\beta}{\delta}=h_2\leq\frac{\alpha}{\eta}.
$$
If $\beta<\eta$, Corollary~\ref{cor:bivariate-lower} gives
$$
D_{(f,X_f)}(h_1,h_2)
\geq
\frac{\eta-\beta}{\delta}
=
\frac{\eta}{\alpha}h_1+h_1-h_2.
$$
Together with the upper bound, this gives the equality.

In the boundary case $\beta=\eta$ for which
$$
h_2=\left(1+\frac{\eta}{\alpha}\right)h_1,
$$
we just need to prove that the corresponding level set is non-empty. For this, it is enough to follow a single branch of the previous construction. Choose $\beta_{m_k}\to\eta$ with $\beta_{m_k}+\varepsilon_k<\eta$ and, at each generation, retain only one $\beta_{m_k}$-good coefficient inside the previous interval. The unique point $x$ in the intersection satisfies
$$
h_f(x)=\frac{\alpha}{\delta},
\qquad
h_{X_f}(x)=\frac{\alpha+\eta}{\delta}.
$$
\end{proof}

\subsection{The randomization effect spectrum}

As a consequence, we obtain the spectrum of the randomization effect. 

\BC\label{cor:gain-spectrum}
Put
$$H_{\max}=\min\left\{\frac{\alpha}{\alpha+\eta},\,\frac{\alpha}{\eta}-\alpha\right\}.$$
Almost surely, for every $H\in[0,H_{\max}]$, one has
$$
D_\Delta(H)=1-\left(1+\frac{\eta}{\alpha}\right)H,
$$
and the level set $\{x:\Delta_f(x)=H\}$ is empty for $H\notin[0,H_{\max}]$. Moreover,
$$
D_\Delta(H)=\sup_{h_2-h_1=H}D_{(f,X_f)}(h_1,h_2) 
 =1- \left( \frac{\eta}{\alpha} +1\right) H.
$$
\EC

\begin{proof}
We work on the event of probability one on which all the preceding results hold. Let
$$E_\Delta(H)=\{x\in[0,1]:h_{X_f}(x)-h_f(x)=H\}.$$
Since $\alpha\leq h_f(x)\leq h_{X_f}(x)\leq\alpha/\eta$ for every $x$, the set $E_\Delta(H)$ is empty if $H<0$ or $H>\alpha/\eta-\alpha$.

Let us first prove the upper bound. The result is clear if $H=0$. Fix $0<H\leq\alpha/\eta-\alpha$ and $0<\varepsilon<H$. If $x\in E_\Delta(H)$, then $h_f(x)\leq\alpha/\eta-H$. Choosing a rational number $q$ such that $h_f(x)\leq q<h_f(x)+\varepsilon$, we obtain
$$
E_\Delta(H)\subseteq\bigcup_{\substack{q\in\mathbb Q\\ \alpha\leq q\leq\alpha/\eta-H+\varepsilon}}\{x:h_f(x)\leq q,\ h_{X_f}(x)\geq q+H-\varepsilon\}.
$$
By Proposition~\ref{prop:general-bivariate-upper} and Corollary~\ref{cor:nu_f}, every set in this union has Hausdorff dimension at most
$$
\frac{\eta}{\alpha}q-H+\varepsilon\leq1-\left(1+\frac{\eta}{\alpha}\right)H+\left(1+\frac{\eta}{\alpha}\right)\varepsilon,
$$
and is empty if this last bound is negative. Since the union is countable, letting $\varepsilon$ tend to zero gives
$$
D_\Delta(H)\leq1-\left(1+\frac{\eta}{\alpha}\right)H
$$
if $H\leq\alpha/(\alpha+\eta)$. If $H>\alpha/(\alpha+\eta)$, we can choose $\varepsilon$ sufficiently small so that the preceding bound is negative. Proposition~\ref{prop:general-bivariate-upper} then shows that every set in the union is empty, hence $E_\Delta(H)=\emptyset$.

For the lower bound if suffices to notice that if $0\leq H\leq H_{\max}$, then by Theorem~\ref{thm:bivariate-spectrum},
$$
D_{\Delta} (H) \geq D_{(f,X_f)}(\frac{\alpha}{\eta}-H,\frac{\alpha}{\eta})= 1-\left(1+\frac{\eta}{\alpha}\right)H.
$$
The conclusion follows.
\end{proof}

As a first consequence, we prove that the set of points where the randomization changes the pointwise H\"older exponent has full Hausdorff dimension.

\begin{proof}[Proof of Theorem~\ref{thm:main_intro}]
Fix $\varepsilon>0$. Choose $H\in(0,H_{\max}]$ sufficiently small so that
$$
1-\left(1+\frac{\eta}{\alpha}\right)H>1-\varepsilon.
$$
By Corollary~\ref{cor:gain-spectrum},
$$
\dim_{\mathcal H}\{x:h_{X_f}(x)\neq h_f(x)\}
\geq D_\Delta(H)
=
1-\left(1+\frac{\eta}{\alpha}\right)H
>1-\varepsilon.
$$
Since $\varepsilon>0$ is arbitrary, we obtain 
$
\dim_{\mathcal H}\{x:h_{X_f}(x)\neq h_f(x)\}=1.$
\end{proof}

As a second consequence, we now show that the bivariate spectrum also allows us to recover the univariate multifractal spectrum of the process $X_f$. Although the pointwise H\"older exponent may be modified on a set of full Hausdorff dimension, the multifractal spectrum remains unchanged.

\BC\label{cor:spectrum_Xf}
Almost surely,
$$
D_{X_f}(h)=\frac{\eta}{\alpha}h
\qquad\text{for every }h\in\left[\alpha,\frac{\alpha}{\eta}\right],
$$
and the level sets $\{x:h_{X_f}(x)=h\}$ are empty outside this interval. In particular, almost surely, for every $h$, one has 
$$
D_{X_f}(h)=D_f(h).
$$
\EC

\begin{proof}
By Corollary~\ref{cor:regularity-increase} and  Proposition \ref{prop:maxregularity}, we know that
$$
\alpha\leq h_{X_f}(x)\leq\frac{\alpha}{\eta}
$$
for every $x\in[0,1]$. Consequently, the level sets of $X_f$ are empty outside $[\alpha,\alpha/\eta]$. Let
$$
h\in\left[\alpha,\frac{\alpha}{\eta}\right].
$$
Taking $h_1=h_2=h$ in Theorem~\ref{thm:bivariate-spectrum}, we obtain
$$
D_{X_f}(h)
\geq D_{(f,X_f)}(h,h)
=
\frac{\eta}{\alpha}h.
$$
On the other hand, Corollary~\ref{cor:nu_f} gives
$$
D_{X_f}(h)\leq\frac{\eta}{\alpha}h.
$$
Hence
$$
D_{X_f}(h)=\frac{\eta}{\alpha}h.
$$
Since Proposition~\ref{prop:det-spectrum} gives the same spectrum for $f$, we conclude that
$$
D_{X_f}(h)=D_f(h)
$$
for every $h$.
\end{proof}

Theorem \ref{thm:RandExp} sums up the main results of 
Section 5.
\BT \label{thm:RandExp}
The process defined by \eqref{eq:X-lacunary} satisfies  almost surely the following properties.
\begin{enumerate}
\item The two marginal multifractal spectra coincide:
\[
 D_f(h)=D_{X_f}(h)=\frac{\eta}{\alpha}h,
 \qquad \forall h\in\left[\alpha,\frac{\alpha}{\eta}\right],
\]
and the level sets are empty outside this interval.
\item One has
\[
 \dim_{\mathcal{H}} \{x:h_f(x)\ne h_{X_f}(x)\}=1.
\]
\item Define the admissible region
\[
 \mathcal R
 =\left\{(h_1,h_2):  \alpha\le h_1\le\frac{\alpha}{\eta},\quad  h_1\le h_2\le\frac{\alpha}{\eta},\quad
 0\le\alpha\left(\frac{h_2}{h_1}-1\right)\leq\eta  \right\}.
\]
Then, for $(h_1,h_2)\in\mathcal R$,
\begin{equation}\label{eq:biv-spectrum-correct}
 D_{(f,X_f)}(h_1,h_2)
 =\left(\frac{\eta}{\alpha} +1\right)h_1-h_2.
\end{equation}
Outside the admissible region, the corresponding level set is empty.
\item For every $H \in [0, \min(\frac{\alpha}{\alpha+\eta},\frac{\alpha}{\eta}-\alpha)] $, one has
$$
D_{\Delta}(H)= 1 - \left( 1+\frac{\eta}{\alpha} \right)H.
$$
\end{enumerate}
\ET

\noindent \textbf{Acknowledgements.} This project was partially supported by the Tournesol program, a Partenariat Hubert Curien (PHC), and by the CNRS International Research Project  FRALEA: Fractalité, Multifractalité et Aléa. C. Esser also acknowledges support from the Fonds Thelam, managed by the King Baudouin Foundation (Belgium). 

\bibliographystyle{plain}
\bibliography{RandomGaussianProcess.bib}

@article{Mandelbrot:74,
     author = {B. Mandelbrot} ,
     title = { Intermittent turbulence in sefl similar cascades},
     journal = {J. Fluid Mech.},
     year = {1974},
     volume={62},
     pages={331--358}
 }

@article{MeyerSellanTaqqu,
  author  = {Meyer, Yves and Sellan, Fabrice and Taqqu, Murad S.},
  title   = {Wavelets, generalized white noise and fractional integration: The synthesis of fractional Brownian motion},
  journal = {Journal of Fourier Analysis and Applications},
  volume  = {5},
  number  = {4},
  pages   = {465--494},
  year    = {1999},
  doi     = {10.1007/BF01261639}
}

@article{pruitt1981growth,
  title={The growth of random walks and L{\'e}vy processes},
  author={Pruitt, William E},
  journal={The Annals of Probability},
  pages={948--956},
  year={1981},
  publisher={JSTOR}
}

@article{AyacheTaqqu,
  author  = {Ayache, Antoine and Taqqu, Murad S.},
  title   = {Multifractional processes with random exponent},
  journal = {Publicacions Matem{\`a}tiques},
  volume  = {49},
  number  = {2},
  pages   = {459--486},
  year    = {2005},
  doi     = {10.5565/PUBLMAT_49205_11}
}

@article{ayache2007wavelet,
  author  = {Ayache, Antoine and Jaffard, St{\'e}phane and Taqqu, Murad S.},
  title   = {Wavelet construction of generalized multifractional processes},
  journal = {Revista Matem{\'a}tica Iberoamericana},
  volume  = {23},
  number  = {1},
  pages   = {327--370},
  year    = {2007},
  doi     = {10.4171/RMI/497}
}

@article{AyacheXiao2005,
  author  = {Ayache, Antoine and Xiao, Yimin},
  title   = {Asymptotic Properties and Hausdorff Dimensions of Fractional Brownian Sheets},
  journal = {Journal of Fourier Analysis and Applications},
  volume  = {11},
  number  = {4},
  pages   = {407--439},
  year    = {2005},
  doi     = {10.1007/s00041-005-4048-3}
}

@article{donoho1994ideal,
  title={Ideal spatial adaptation by wavelet shrinkage},
  author={Donoho, David L. and Johnstone, Iain M.},
  journal={Biometrika},
  volume={81},
  number={3},
  pages={425--455},
  year={1994}
}

@article{donoho1995wavelet,
  title={Wavelet shrinkage: Asymptopia?},
  author={Donoho, David L. and Johnstone, Iain M. and Kerkyacharian, G{\'e}rard and Picard, Dominique},
  journal={Journal of the Royal Statistical Society: Series B (Methodological)},
  volume={57},
  number={2},
  pages={301--337},
  year={1995}
}

@article{donoho1995denoising,
  title={De-noising by soft-thresholding},
  author={Donoho, David L.},
  journal={IEEE Transactions on Information Theory},
  volume={41},
  number={3},
  pages={613--627},
  year={1995}
}

@article{donoho1995adapting,
  title={Adapting to unknown smoothness via wavelet shrinkage},
  author={Donoho, David L. and Johnstone, Iain M.},
  journal={Journal of the American Statistical Association},
  volume={90},
  number={432},
  pages={1200--1224},
  year={1995}
}

@book{talagrand2014upper,
  title={Upper and lower bounds for stochastic processes},
  author={Talagrand, Michel},
  year={2014},
  publisher={Springer}
}

@book{ledoux1991probability,
  title={Probability in Banach Spaces: isoperimetry and processes},
  author={Ledoux, Michel and Talagrand, Michel},
  volume={23},
  year={1991},
  publisher={Springer}
}

@article{grebenkov2016multiscale,
  title={A multiscale guide to Brownian motion},
  author={Grebenkov, Denis S and Belyaev, Dmitry and Jones, Peter W},
  journal={Journal of Physics A: Mathematical and Theoretical},
  volume={49},
  number={4},
  pages={043001},
  year={2016},
  publisher={IOP Publishing}
}

@article{ciesielski1961holder,
  title={H{\"o}lder conditions for realizations of Gaussian processes},
  author={Ciesielski, Zbigniew},
  journal={Transactions of the American Mathematical Society},
  pages={403--413},
  year={1961},
  publisher={JSTOR}
}

@article{orey1974often,
  title={How often on a Brownian path does the law of iterated logarithm fail?},
  author={Orey, Steven and Taylor, S James},
  journal={Proceedings of the London Mathematical Society},
  volume={3},
  number={1},
  pages={174--192},
  year={1974},
  publisher={Wiley Online Library}
}

@incollection{khoshnevisan2007fast,
  title={Fast sets and points for fractional Brownian motion},
  author={Khoshnevisan, Davar and Shi, Zhan},
  booktitle={S{\'e}minaire de Probabilit{\'e}s XXXIV},
  pages={393--416},
  year={2007},
  publisher={Springer}
}

@incollection{fernique2006regularite,
  title={Regularit{\'e} des trajectoires des fonctions al{\'e}atoires gaussiennes},
  author={Fernique, Xavier},
  booktitle={Ecole d’Et{\'e} de Probabilit{\'e}s de Saint-Flour IV—1974},
  pages={1--96},
  year={2006},
  publisher={Springer}
}

@incollection{dudley2010sample,
  title={Sample functions of the Gaussian process},
  author={Dudley, Richard M},
  booktitle={Selected works of RM Dudley},
  pages={187--224},
  year={2010},
  publisher={Springer}
}

@inproceedings{jain1974sufficient,
  title={Sufficient conditions for the continuity of stationary Gaussian processes and applications to random series of functions},
  author={Jain, Naresh C and Marcus, Michael B},
  publisher={Annales de l'institut Fourier},
  volume={24(2)},
  pages={117--141},
  year={1974}
}

@book{marcus2016random,
  title={Random Fourier Series with Applications to Harmonic Analysis.(AM-101)},
  series={Annals of Mathematic Studies},
  author={Marcus, Michael B and Pisier, Gilles},
  year={2016},
  publisher={Princeton University Press}
}

@ARTICLE{Beresn:06,
  author = {V. Beresnevich and S. Velani},
  title = {A mass transference principle and the {D}uffin-{S}chaeffer conjecture
	for {H}ausdorff measures},
  journal = {Ann. of Math.},
  year = {2006},
  volume = {164},
  pages = {971--992}
}

@ARTICLE{Jaffard:02,
  author = {J.M. Aubry and S. Jaffard},
  title = {Random wavelet series},
  journal = {Comm. Math. Phys.},
  year = 2002,
  volume = 227,
  pages = {483--514}
}

@article{benassi1997elliptic,
  title={Elliptic Gaussian random processes},
  author={Benassi, Albert and Jaffard, St{\'e}phane and Roux, Daniel},
  journal={Revista matem{\'a}tica iberoamericana},
  volume={13},
  number={1},
  pages={19--90},
  year={1997}
}

@article{horst2024besov,
  title={Besov regularity of random wavelet series},
  author={Horst, Andreas and Jahn, Thomas and Voigtlaender, Felix},
  journal={arXiv preprint arXiv:2411.18155},
  year={2024}
}

@article{esser2025p,
  title={The p-spectrum of Random Wavelet Series},
  author={Esser, C{\'e}line and Lambert, Thelma and B{\'e}atrice, Vedel},
  journal={arXiv preprint arXiv:2510.00622},
  year={2025}
}

@ARTICLE{Bastin:16,
  author = {F. Bastin and C. Esser and S. Jaffard},
  title = {Large deviation spectra based on wavelet leaders},
  journal = {Rev. Matem. Iberoamer.},
  year={2016},
  volume = {32 (3)},
  pages = {859--890}
}

@ARTICLE{Esser:17,
author = {C. Esser and T. Kleyntssens and S. Nicolay},
year = {2017},
pages = {269--291},
title = {A multifractal formalim for non-concave and non-increasing spectra: The leaders profile method},
journal = {Appl. Comput. Harmon. Anal.},
volume = {43 (2)}
}

@BOOK{Falconer:86,
  title = {The Geometry of Fractal Sets},
  publisher = {Cambridge University Press},
  year = {1986},
  author = {K. Falconer},
  isbn = {9780521337052},
}

@article{ayache2021monofractality,
  title={On the monofractality of many stationary continuous gaussian fields},
  author={Ayache, Antoine},
  journal={Journal of Functional Analysis},
  volume={281},
  number={7},
  pages={109111},
  year={2021},
  publisher={Elsevier}
}

@article{ayache2013continuous,
  title={Continuous Gaussian multifractional processes with random pointwise H{\"o}lder regularity},
  author={Ayache, Antoine},
  journal={Journal of Theoretical Probability},
  volume={26},
  number={1},
  pages={72--93},
  year={2013},
  publisher={Springer}
}

@ARTICLE{Jaffard:04b,
  author = {S. Jaffard},
  title = {Wavelet techniques in multifractal analysis, Fractal Geometry and
	Applications: A Jubilee of {B}enoit {M}andelbrot},
  journal = {Proceedings of Symposia in Pure Mathematics},
  year = {2004},
  volume = {72},
  pages = {91--151},
  editor = {M. Lapidus and M. van Frankenhuijsen}
}

@ARTICLE{Jaffard:99,
  author = {S Jaffard},
  title = {The multifractal nature of {L}{\'e}vy processes},
  journal = {Probab. Theory Relat. Fields},
  year = {1999},
  volume = {114},
  pages = {207--227}
}

@ARTICLE{Jaffard:00b,
author={S. Jaffard}, 
title={Lacunary wavelet series},
journal={Annals of Applied Probability.},
year={2000},
volume= {Vol. 10, No. 1}, 
pages={pp. 313-329},
publisher ={Institute of Mathematicals Society}
}

@article{kaufman1975large,
  title={Large increments of Brownian motion},
  author={Kaufman, R},
  journal={Nagoya Mathematical Journal},
  volume={56},
  pages={139--145},
  year={1975},
  publisher={Cambridge University Press}
}

@article{kahane1974irregularite,
  title={Sur l’irr{\'e}gularit{\'e} locale du mouvement brownien},
  author={Kahane, Jean-Pierre},
  journal={CR Acad. Sci. Paris},
  volume={278},
  pages={331--333},
  year={1974}
}

@article{LM86,
Author={P.G. Lemari\'e and Y. Meyer},
title = {Ondelettes et bases hilbertiennes},
journal={Rev. Mat. Iberoamericana},
year={1986},
volume={1}
}

@article{jaffard1996local,
  title={Local regularity of nonsmooth wavelet expansions and application to the Polya function},
  author={Jaffard, St{\'e}phane and Mandelbrot, Benoit B},
  journal={advances in mathematics},
  volume={120},
  number={2},
  pages={265--282},
  year={1996},
  publisher={Elsevier}
}

@BOOK{Meyer:95,
  title = {Ondelettes et {O}p\'erateurs},
  publisher = {Hermann, Paris},
  year = {1990},
  author = {Y. Meyer},
  volume = {1}
}

@article{stromberg2006modified,
  title={A modified Franklin system and higher-order spline systems on Rn as unconditional bases for Hardy spaces},
  author={Str{\"o}mberg, Jan-Olov},
  journal={Fundamental Paper in Wavelet Theory},
  pages={197--215},
  year={2006},
  publisher={Princeton University Press}
}

@article {SaesSeuret,
    AUTHOR = {G. Sa\"{e}s and S. Seuret},
     TITLE = {Multifractal analysis of sums of random pulses},
   JOURNAL = {Math. Proc. Cambridge Philos. Soc.},
  FJOURNAL = {Mathematical Proceedings of the Cambridge Philosophical
              Society},
    VOLUME = {175},
      YEAR = {2023},
    NUMBER = {3},
     PAGES = {569--593},
      ISSN = {0305-0041,1469-8064},
   MRCLASS = {28A80 (26A15 28A78 60G17 60L30)},
  MRNUMBER = {4655525},
       DOI = {10.1017/s0305004123000324},
       URL = {https://doi.org/10.1017/s0305004123000324},
}

@article{EL22,
author={Esser, C. and Loosveldt, L.},
title={{S}low, ordinary and rapid points for {G}aussian {W}avelets {S}eries and application to {F}ractional {B}rownian {M}otions},
journal={Latin American Journal of Probability and Mathematical Statistics},
volume={19},
pages={1479-1495},
year={2022}
}

@BOOK{Dau92,
	AUTHOR={I. Daubechies},
	TITLE={Ten Lectures on Wavelets},
	PUBLISHER={Society for Industrial and Applied Mathematics},
	YEAR={1992},
	ADDRESS={Philadelphia, PA}
}

@book{jaffard1996wavelet,
  title={Wavelet methods for pointwise regularity and local oscillations of functions},
  author={Jaffard, St{\'e}phane and Meyer, Yves},
  volume={587},
  year={1996},
  publisher={American Mathematical Soc.}
}

@BOOK{Mey90I,
	AUTHOR={Y. Meyer},
	TITLE={Ondelettes et Op\'erateurs},
	PUBLISHER={Hermann},
	YEAR={1990},
	ADDRESS={Paris},
	NOTE={English translation, {\it Wavelets and operators}, Cambridge University Press, 1992.}
}

@BOOK{Mallat1998,
  title = {A Wavelet Tour of Signal Processing},
  publisher = {Academic Press},
  year = {1998},
  author = {Mallat, S.},
  address = {San Diego, CA}
}

@article{EJV24,
  title={Regularity properties of random wavelet series},
  author={Esser, C{\'e}line and Jaffard, St{\'e}phane and Vedel, B{\'e}atrice},
  journal={Theory of Probability and Mathematical Statistics},
  volume={110},
  pages={31--53},
  year={2024}
}

\end{document}